# Generation of Curvilinear Coordinates

Dedicated to Professor Emeritus Yoshiyuki Yamamoto, The University of Tokyo

Hiroshi Isshiki, Institute of Mathematical Analysis, Osaka, Japan,
isshiki@dab.hi-ho.ne.jp

Daisuke Kitazawa, Institute of Industrial Science,
The University of Tokyo, dkita@iis.u-tokyo.ac.jp

**Abstract**

The authors have discussed the method and merit of introducing the curvilinear coordinates into numerical analysis and have shown some numerical examples. However, they used analytical functions for the mappings in the examples and didn't mention how to generate a mapping numerically between the physical and mapped coordinates. There is an analytical method using the solution of Dirichlet problem of Poisson equation. The authors present an algebraic method using interpolation of mapping function values at discrete points based on the least square method. Not only the mapping of two coordinates but also the application to the solution of boundary value problem is given.

## 1. Introduction

Curvilinear coordinates might be more convenient in many cases than Cartesian coordinates. In order to treat rapid change of the solution function, we usually use finer mesh to follow the change. Curvilinear coordinates would be able to handle fixed discontinuity, moving discontinuity in more natural way. Namely, a curvilinear coordinates could generate fine mesh in the neighborhood of fixed and moving discontinuities and increase the accuracy of the solution.

Furthermore, if we introduce curvilinear coordinates, we can transform a non-square region into square one and can introduce a regular mesh in the transformed coordinates. Furthermore, in numerical calculation, a curved boundary is usually approximated by a jagged or non-smooth boundary. We could treat the boundary more reasonably.

The authors have discussed the method and merit of introducing the curvilinear coordinates into numerical analysis and have shown some numerical examples [1]. However, they used analytical functions for the mappings in the examples and didn't mention how to generate a mapping numerically between the physical and mapped coordinates.

There is an analytical method using the solution of Dirichlet problem of Poisson equation [2, 3]. In the present paper, the authors present an algebraic method using interpolation of mapping function values at discrete points based on the Least Square Method (LSM).

The analytical method using Dirichlet problem of Poisson equation could be summarized as follows. The merit of this method is that discrete function values are specified only on the boundary of the region. This would be very convenient for the 3D mesh generation. However, the curvilinear coordinates are given as a set of discrete values, and the concentration of the mesh points are indirectly specified by the inhomogeneous term of Poisson equation.

On the other hand, the algebraic method interpolates discretely given function values using LSM. The merit of this method is that the curvilinear coordinates are given as continuous functions. This property is very convenient for calculating derivatives of the mapping functions. If the mapping data are specified using CAD, the knowhow of the specialist could be applied to the mapping. However, the discrete function values must be specified not only on the boundary but also in the interior region

Not only the mapping of two coordinates but also the application to the solution of boundary value problem is given.

## 2. Mapping from a non-rectangular region to a rectangular one

Generally speaking, there are two types o coordinates. One is Cartesian type based on “Cartesian coordinates” as shown in Figs. 1(a) and 2(a). The other is Polar type based on “Polar coordinates” as shown in Figs. 1(b) and 2(b). Let a region in Cartesian coordinates $(x, y)$ be mapped into a rectangular region in curvilinear coordinates $(\xi,\eta)$ :

$$\xi = \xi(x, y)\,,\quad \eta = \eta(x, y) \quad \text{or} \quad x = x(\xi,\eta)\,,\quad y = y(\xi,\eta)\,. \tag{1}$$

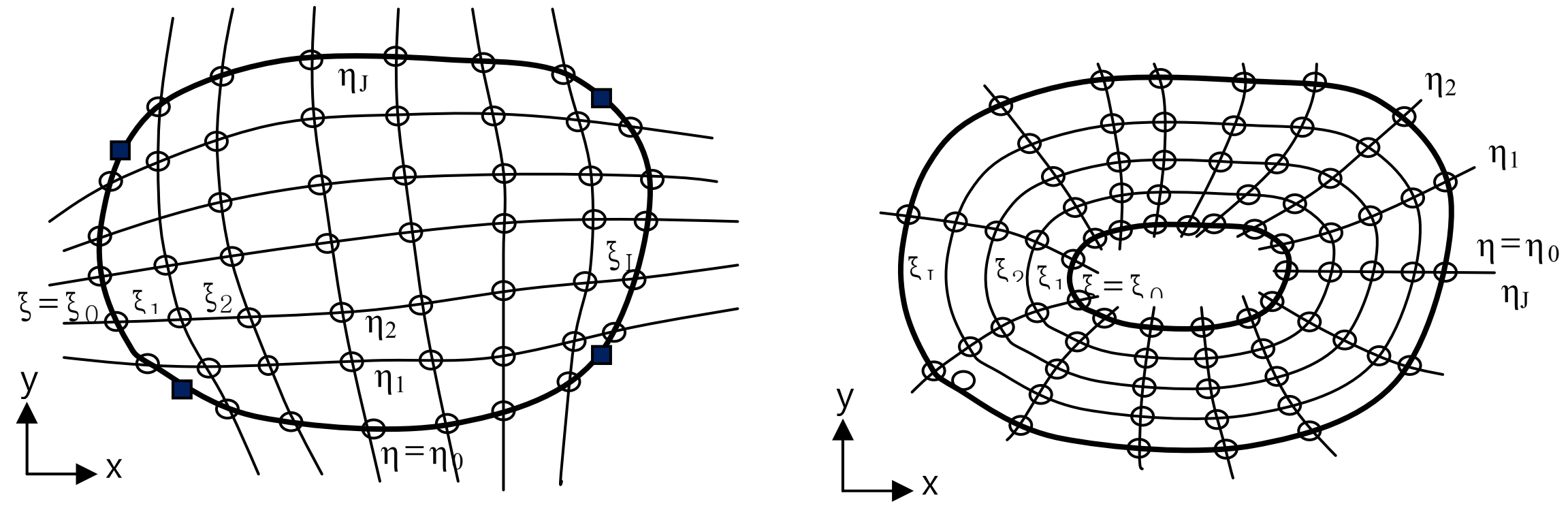

(a) Cartesian type in simply connected region (b) Polar type in doubly connected region

Fig. 1. Typical curvilinear coordinates.

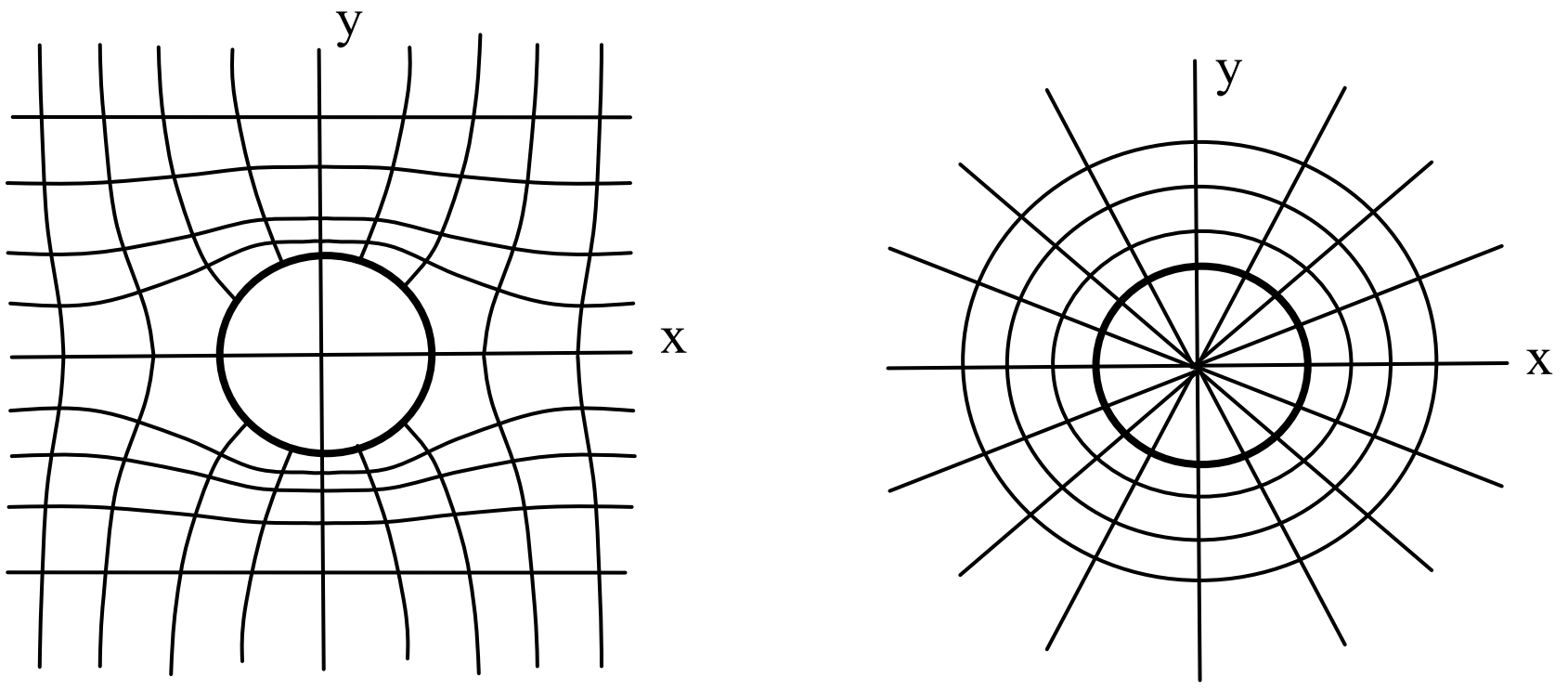

(a) Cartesian type (b) Polar type

Fig. 2. Two types of curvilinear coordinates in doubly connected region.

For the normal transformation from the physical plane into the mapped plane, we assume

$$\xi = \sum_{m=0}^{M}\sum_{n=0}^{m+1} A_{mn} x^{m-n} y^{n}\,,\quad \eta = \sum_{m=0}^{M}\sum_{n=0}^{m+1} B_{mn} x^{m-n} y^{n}\,. \tag{2a,b}$$

We define mesh points in the mapped region $(\xi_i, \eta_j)$:

$$\xi_i = \xi_0 + \frac{\xi_I - \xi_0}{I} i\,, \quad i = 0,1,\cdots,I\ ; \qquad \eta_i = \eta_0 + \frac{\eta_J - \eta_0}{J} j\,, \quad j = 0,1,\cdots,J \tag{3a}$$

and the corresponding mesh points $(x_{ij}, y_{ij})$:

$$x_{ij} = x(\xi_i, \eta_j)\,, \quad y_{ij} = y(\xi_i, \eta_j)\,. \tag{3b}$$

Then, we obtain for $i = 0,1,\cdots,I$:

$$\xi_i = \sum_{m=0}^{M} \sum_{n=0}^{m+1} A_{mn} x_{ij}^{\,m-n} y_{ij}^{\,n}\,, \quad j = 0,1,\cdots,J\,. \tag{4a}$$

Equation (4a) can be written using matrix notation for $i = 0,1,\cdots,I$

$$\begin{bmatrix} \xi_i \\ \xi_i \\ \cdots \\ \xi_i \end{bmatrix} = \begin{bmatrix} 1 & x_{i0} & y_{i0} & x_{i0}^{\,2} & x_{i0}y_{i0} & y_{i0}^{\,2} & \cdots \\ 1 & x_{i1} & y_{i1} & x_{i1}^{\,2} & x_{i1}y_{i1} & y_{i1}^{\,2} & \cdots \\ \cdots & \cdots & \cdots & \cdots & \cdots & \cdots & \cdots \\ 1 & x_{iJ} & y_{iJ} & x_{iJ}^{\,2} & x_{iJ}y_{iJ} & y_{iJ}^{\,2} & \cdots \end{bmatrix} \begin{bmatrix} A_{00} \\ A_{10} \\ A_{01} \\ A_{20} \\ A_{11} \\ A_{02} \\ \vdots \end{bmatrix}. \tag{4b}$$

Equation (4) or ((5) is solved using Least Square Method (LSM). The total number of mesh points $(I+1)\times(J+1)$ must be equal to or greater than the total number of the unknowns $(M+1)(M+2)/2$.

Similarly, we have for $j = 0,1,\cdots,J$,

$$\eta_j = \sum_{m=0}^{M} \sum_{n=0}^{m+1} B_{mn} x_{ij}^{\,m-n} y_{ij}^{\,n}\,, \quad i = 0,1,\cdots,I\,. \tag{5a}$$

Using matrix notation, we have for $j = 0,1,\cdots,J$

$$\begin{bmatrix} \eta_j \\ \eta_j \\ \cdots \\ \eta_j \end{bmatrix} = \begin{bmatrix} 1 & x_{0j} & y_{0j} & x_{0j}^{\,2} & x_{0j}y_{0j} & y_{0j}^{\,2} & \cdots \\ 1 & x_{1j} & y_{1j} & x_{1j}^{\,2} & x_{1j}y_{1j} & y_{1j}^{\,2} & \cdots \\ \cdots & \cdots & \cdots & \cdots & \cdots & \cdots & \cdots \\ 1 & x_{Ij} & y_{Ij} & x_{Ij}^{\,2} & x_{Ij}y_{Ij} & y_{Ij}^{\,2} & \cdots \end{bmatrix} \begin{bmatrix} B_{00} \\ B_{10} \\ B_{01} \\ B_{20} \\ B_{11} \\ B_{02} \\ \vdots \end{bmatrix}. \tag{5b}$$

Equation (4) or (5) is solved using Least Square Method (LSM). The total number of mesh points $(I+1)\times(J+1)$ must be equal to or greater than the total number of the unknowns $(M+1)(M+2)/2$.

## 3. Inverse transformation from a mapped region to a physical one

For the inverse transformation from the mapped plane into the physiical plane, we assume

$$x = \sum_{m=0}^{M} \sum_{n=0}^{m+1} A'_{mn} \xi^{m-n} \eta^{n}\,, \quad y = \sum_{m=0}^{M} \sum_{n=0}^{m+1} B'_{mn} \xi^{m-n} \eta^{n}\,. \tag{6a,b}$$

Then, we obtain for $i = 0,1,\cdots,I$:

$$x_{ij} = \sum_{m=0}^{M} \sum_{n=0}^{m+1} A'_{mn} \xi_i^{\,m-n} \eta_j^{\,n}\,, \quad j = 0,1,\cdots,J \tag{7a}$$

Equation (9) can be written using matrix notation for $i = 0,1,\cdots,I$

$$\begin{bmatrix} x_{i0} \\ x_{i1} \\ \cdots \\ x_{iJ} \end{bmatrix} = \begin{bmatrix} 1 & \xi_i & \eta_0 & \xi_i^2 & \xi_i\eta_0 & \eta_0^2 & \cdots \\ 1 & \xi_i & \eta_1 & \xi_i^2 & \xi_i\eta_1 & \eta_1^2 & \cdots \\ \cdots & \cdots & \cdots & \cdots & \cdots & \cdots & \cdots \\ 1 & \xi_i & \eta_J & \xi_i^2 & \xi_i\eta_J & \eta_J^2 & \cdots \end{bmatrix} \begin{bmatrix} A'_{00} \\ A'_{10} \\ A'_{01} \\ A'_{20} \\ A'_{11} \\ A'_{02} \\ \vdots \end{bmatrix}. \tag{7b}$$

Similarly, we have for $j = 0,1,\cdots,J$,

$$y_{ij} = \sum_{m=0}^{M}\sum_{n=0}^{m+1} B'_{mn}\xi_i^{m-n}\eta_j^{n}, \quad i = 0,1,\cdots,I. \tag{8a}$$

Using matrix notation, we have for $j = 0,1,\cdots,J$

$$\begin{bmatrix} y_{0j} \\ y_{1j} \\ \cdots \\ y_{Ij} \end{bmatrix} = \begin{bmatrix} 1 & \xi_0 & \eta_j & \xi_0^2 & \xi_0\eta_j & \eta_j^2 & \cdots \\ 1 & \xi_1 & \eta_j & \xi_1^2 & \xi_1\eta_j & \eta_j^2 & \cdots \\ \cdots & \cdots & \cdots & \cdots & \cdots & \cdots & \cdots \\ 1 & \xi_I & \eta_j & \xi_I^2 & \xi_I\eta_j & \eta_j^2 & \cdots \end{bmatrix} \begin{bmatrix} B'_{00} \\ B'_{10} \\ B'_{01} \\ B'_{20} \\ B'_{11} \\ B'_{02} \\ \vdots \end{bmatrix}. \tag{8b}$$

## 4. Numerical examples

### 4.1. An example of Cartesian type mapping

We consider a simple example. The mesh points are given in Table 1 and shown in Fig. 3.

Table 1. Given mesh points ($I$=4, $J$=4, $M$=3).

| $\xi$ | $\eta$ | $x$ | $y$ | $\xi$ | $\eta$ | $x$ | $y$ |
|---|---|---|---|---|---|---|---|
| 0 | 0 | 1 | 2 | 0.75 | 0 | 6.85 | 1.25 |
| 0 | 0.25 | 1.5 | 3.35 | 0.75 | 0.25 | 6.8 | 3.2 |
| 0 | 0.5 | 2 | 4.5 | 0.75 | 0.5 | 6.75 | 4.85 |
| 0 | 0.27 | 2.5 | 5.75 | 0.75 | 0.27 | 6.7 | 6.25 |
| 0 | 1 | 3 | 7 | 0.75 | 1 | 6.6 | 7.75 |
| 0.25 | 0 | 3.1 | 1.7 | 1 | 0 | 9 | 1 |
| 0.25 | 0.25 | 3.35 | 3.2 | 1 | 0.25 | 8.7 | 3.2 |
| 0.25 | 0.5 | 3.75 | 4.6 | 1 | 0.5 | 8.4 | 5 |
| 0.25 | 0.27 | 4 | 5.9 | 1 | 0.27 | 8.2 | 6.4 |
| 0.25 | 1 | 4.3 | 7.25 | 1 | 1 | 8 | 8 |
| 0.5 | 0 | 5 | 1.4 | | | | |
| 0.5 | 0.25 | 5.2 | 3.2 | | | | |
| 0.5 | 0.5 | 5.3 | 4.75 | | | | |
| 0.5 | 0.27 | 5.3 | 6.1 | | | | |
| 0.5 | 1 | 5.5 | 7.5 | | | | |

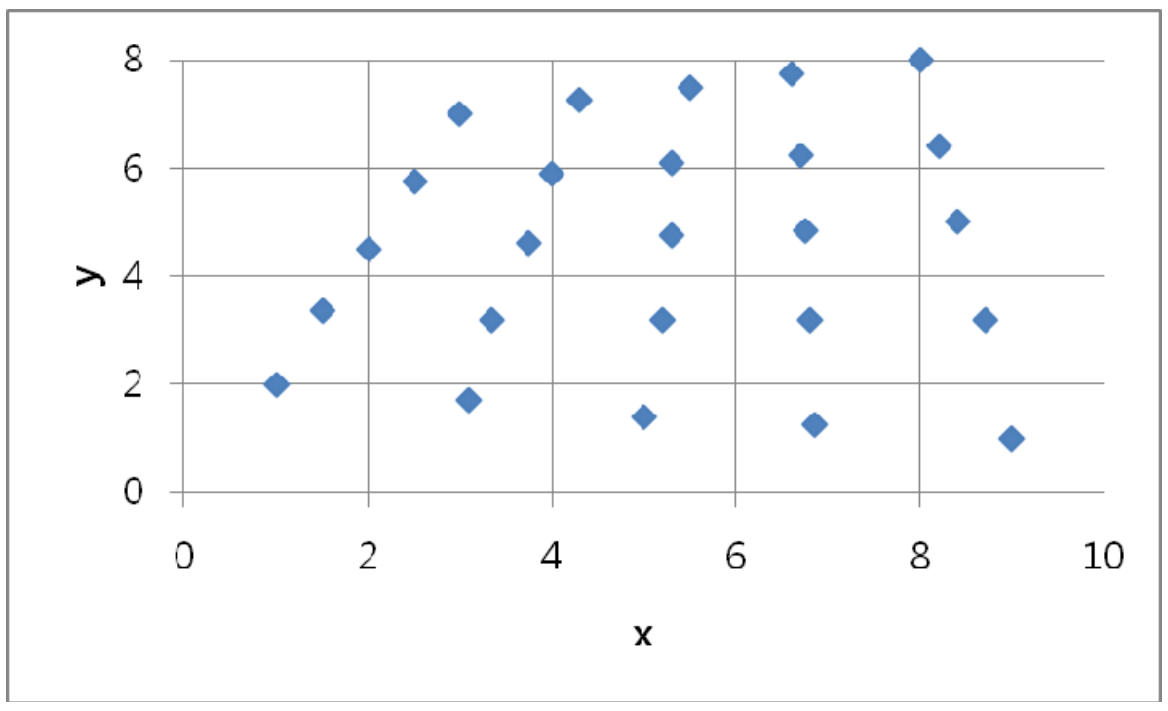

Fig. 3. Physical plane.

In the calculation below, the following parameters are used:

$$I = 4,\quad J = 4,\quad M = 3. \tag{9}$$

The solution of the normal mapping $A_{mn}$ and $B_{mn}$ are determined by Eqs. (4) and (5) and are shown in Table 2. The mapping results are shown in Fig. 4. Small amounts of errors are recognized in Fig. 4(b).

Table 2. $A_{mn}$ and $B_{mn}$ ($I$=4, $J$=4, $M$=3).

| $m$ | $n$ | $A_{mn}$ | $B_{mn}$ |
|---|---|---|---|
| 0 | 0 | 0.012661 | -0.35625 |
| 1 | 0 | 0.042116 | 0.070531 |
| 1 | 1 | -0.02541 | 0.122901 |
| 2 | 0 | 0.016011 | -0.00651 |
| 2 | 1 | 0.006537 | -0.01321 |
| 2 | 2 | -0.00659 | 0.018831 |
| 3 | 0 | -0.00094 | 0.000247 |
| 3 | 1 | -0.00047 | 0.000597 |
| 3 | 2 | 0.001222 | -0.00017 |
| 3 | 3 | -0.00012 | -0.00094 |

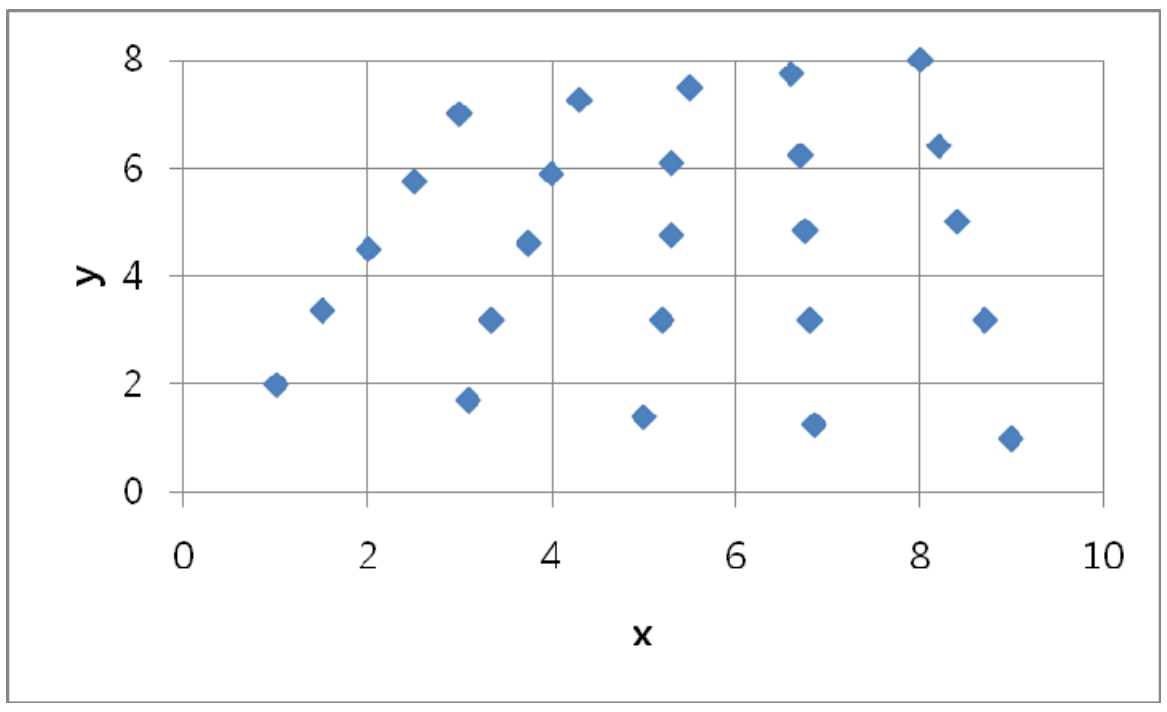

(a) Physical plane $(x,y)$

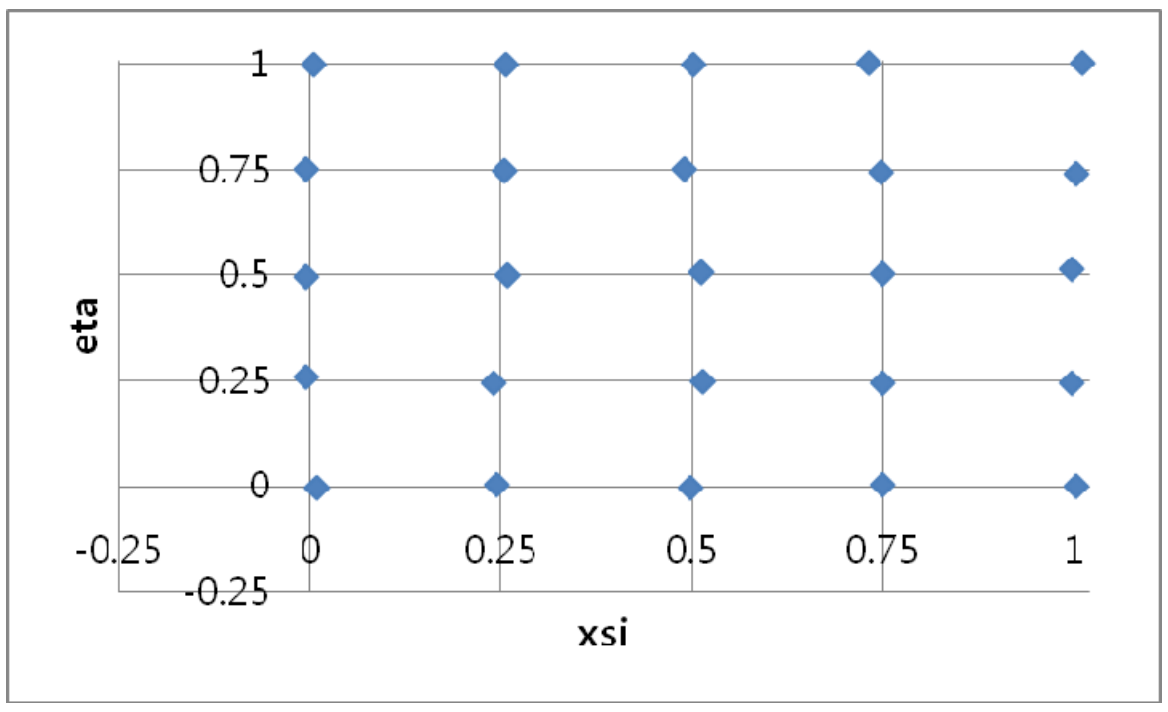

(b) Mapped plane $(\xi,\eta)$

Fig. 4. Normal mapping $\xi=\xi(x,y)$, $\eta=\eta(x,y)$.

The solution of the inverse mapping $A'_{mn}$ and $B'_{mn}$ are determined by Eqs. (7) and (8) and are shown in Table 2. The mapping results are shown in Fig. 5. The accuracy of the inverse mapping is high. The generalization ability is also satisfactory as shown in Fig. 6.

Table 3. $A'_{mn}$ and $B'_{mn}$ ($I$=4, $J$=4, $M$=3).

| $m$ | $n$ | $A'_{mn}$ | $B'_{mn}$ |
|---|---|---|---|
| 0 | 0 | 0.981 | 1.988714 |
| 1 | 0 | 9.122476 | -1.44229 |
| 1 | 1 | 2.196286 | 6.027714 |
| 2 | 0 | -3.16571 | 0.8 |
| 2 | 1 | -3.39543 | 3.998857 |
| 2 | 2 | -0.35429 | -2.75429 |
| 3 | 0 | 2.026667 | -0.32 |
| 3 | 1 | 0.16 | -0.38857 |
| 3 | 2 | 0.251429 | -1.71429 |
| 3 | 3 | 0.16 | 1.76 |

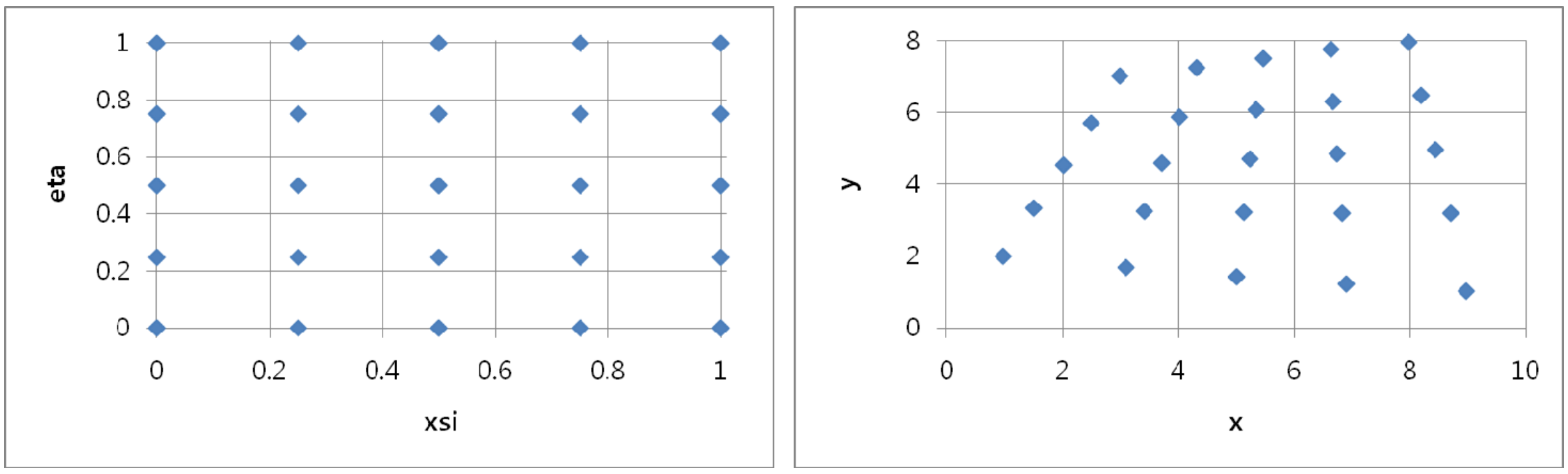

(a) Mapped plane $(\xi,\eta)$ (b) Physical plane $(x,y)$

Fig. 5. Inverse mapping $x=x(\xi,\eta)$, $y=y(\xi,\eta)$.

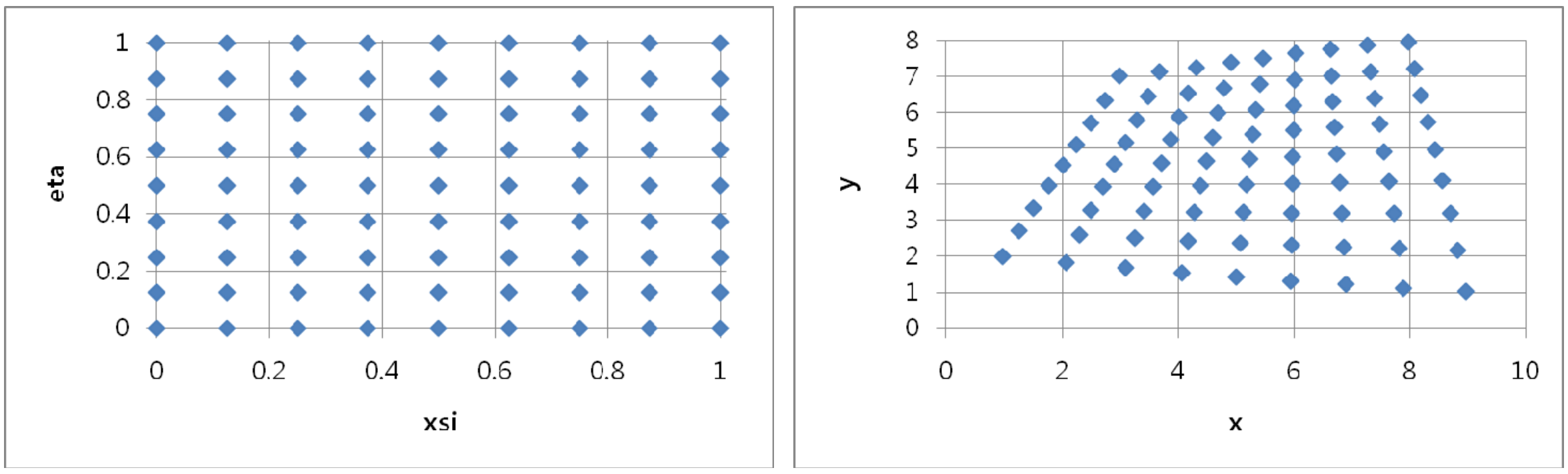

(a) Mapped plane $(\xi,\eta)$ (b) Physical plane $(x,y)$

Fig. 6. Generalization ability of inverse mapping from mapped plane to physical plane ($I$=4, $J$=4, $M$=3).

If we use $M = 5$, the accuracy of Figs. 4(b) and 5(b) are increased, but that of Fig. 6(b) are decreased as shown in Fig. 7. This means that too big $M$ is not good. The proper selection of $M$ seems important.

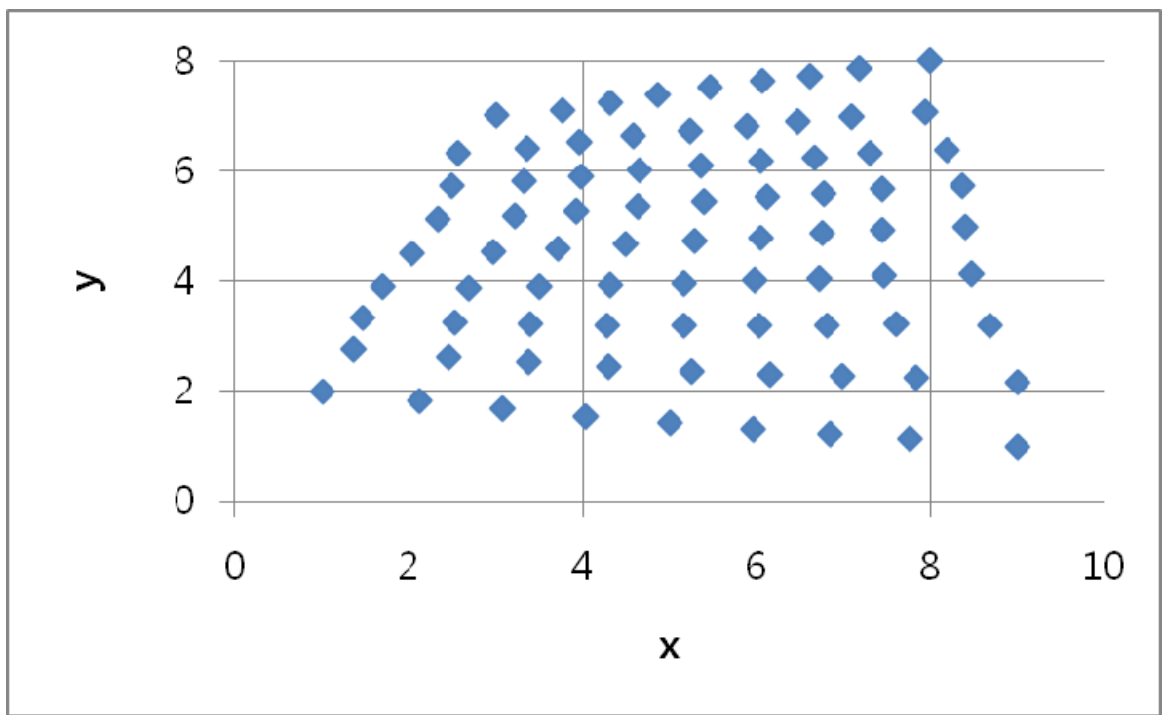

Fig. 7. Mapping results when $M$=5.

If we increase mesh points as shown in Fig. 8, we obtain better result as shown in Figs. 9 and 10.

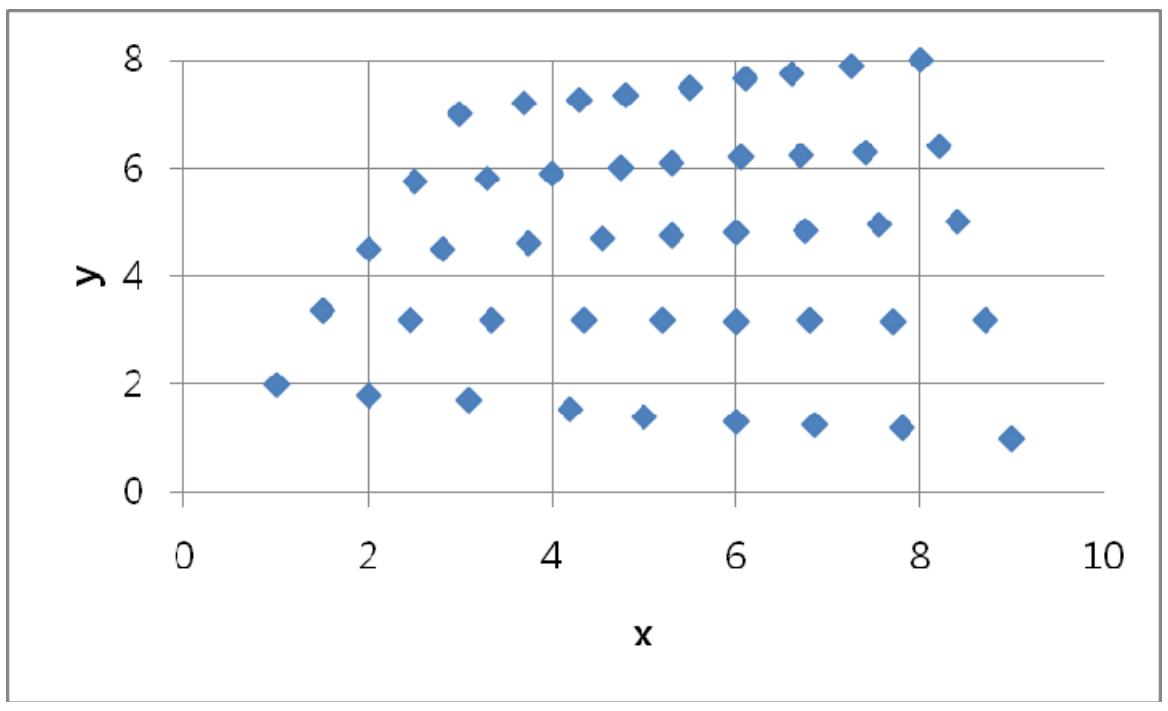

Fig. 8. Physical plane ($I$=8, $J$=4).

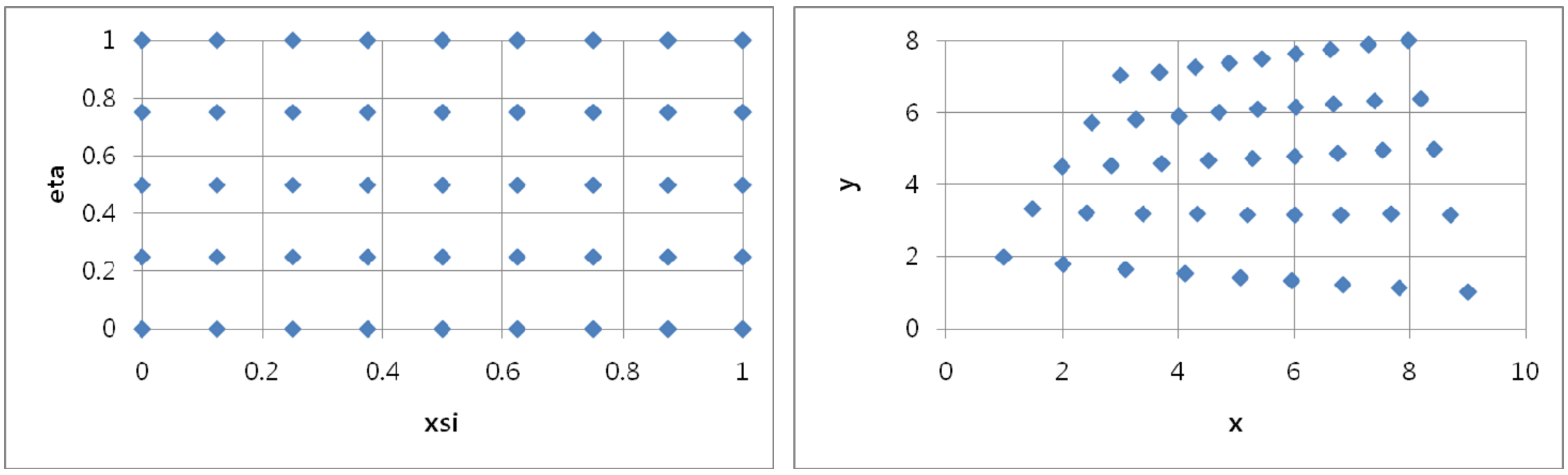

(a) Mapped plane ($\xi$,$\eta$) (b) Physical plane ($x$,$y$)

Fig. 9. Inverse mapping $x$=$x$($\xi$,$\eta$), $y$=$y$($\xi$,$\eta$) ($I$=8, $J$=4, $M$=5).

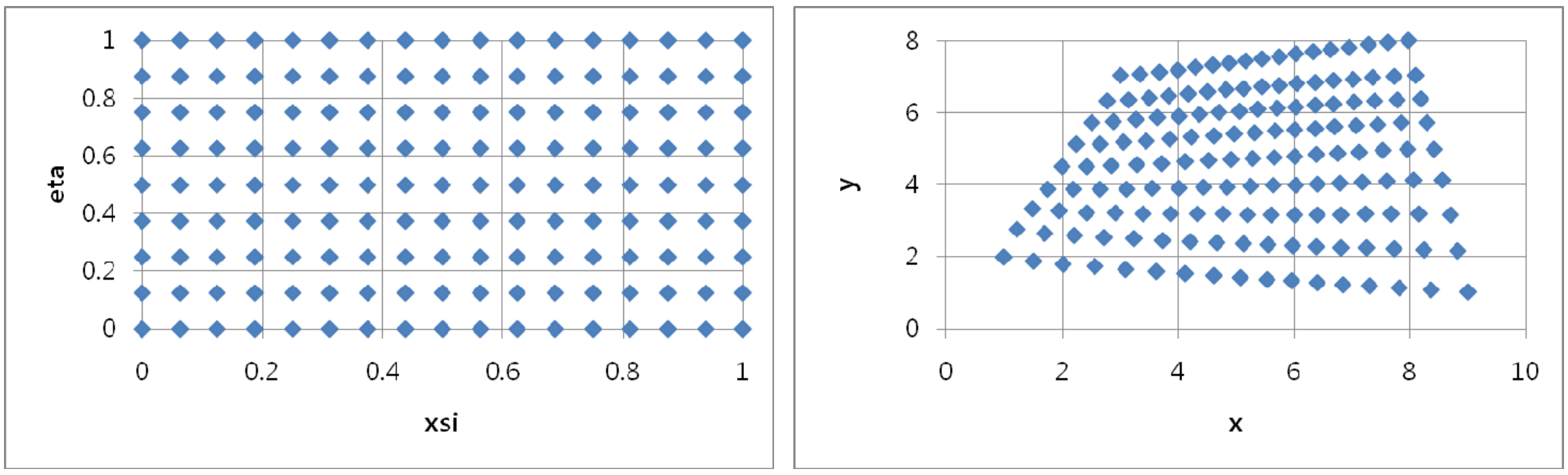

(a) Mapped plane ($\xi$,$\eta$) (b) Physical plane ($x$,$y$)

Fig. 10. Generalization ability of inverse mapping from mapped plane to physical plane ($I$=8, $J$=4, $M$=5).

### 4.2. An example of polar type mapping

(1) Mapping of a full circle

Given points for the rectangular coordinates in the physical plane are specified as

$$x_{ij} = \xi_i \cos \eta_j ,\quad y_{ij} = \xi_i \sin \eta_j ,\quad i = 0,1,\cdots,I ,\quad j = 0,1,\cdots,J . \tag{10a}$$

And given points for the polar coordinates in the mapped plane are specified as

$$\xi_i = 1 + \frac{2-1}{I} i ,\quad i = 0,1,\cdots,I ;\quad \eta_j = 0 + \frac{2\pi - 0}{J} j ,\quad j = 0,1,\cdots,J . \tag{10b}$$

In the polar type mapping, we must notice that the line with $j = 0$ or $\eta_0 = 0$ and the line with $j = J$ or $\eta_J = 2\pi$ coincide and corresponds to $y = 0$, $a < x < R$. This means that two values are specified on $y = 0$, $a < x < R$. Hence, in normal mapping or mapping from $(x, y)$ to $(\xi, \eta)$, $\eta_j$ must be modified as

$$\eta^*_j = \eta_j - \frac{2\pi}{J} j, \tag{11a}$$

and, after the mapping, $\eta_j$ is obtained as

$$\eta_j = \eta^*_j + \frac{2\pi}{J} j. \tag{11b}$$

In the calculations below, the following parameters are used:

$$I = 4, \quad J = 16, \quad M = 7. \tag{12}$$

The solution of the normal mapping $A_{mn}$ and $B_{mn}$ are determined by Eqs. (4) and (5). The mapping results are shown in Figs. 11 and 12. The accuracy of the normal mapping is high as shown in Fig. 13.

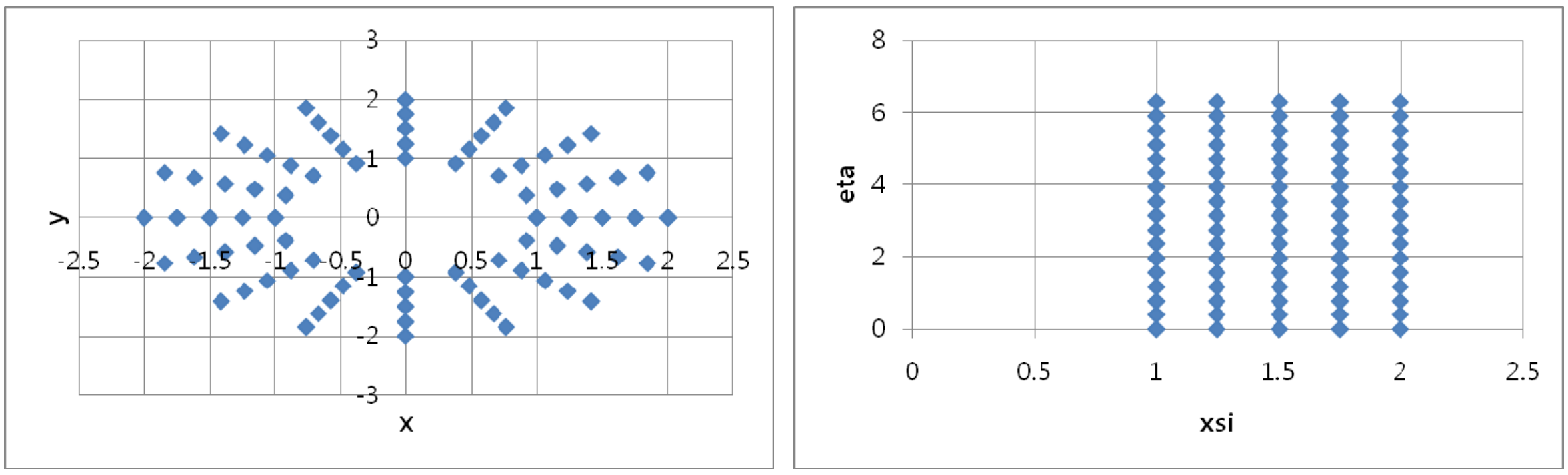

(a) Physical plane $(x,y)$ (b) Mapped plane $(\xi,\eta)$

Fig. 11. Normal mapping $\xi=\xi(x,y)$, $\eta=\eta(x,y)$ ($I$=4, $J$=16, $M$=7).

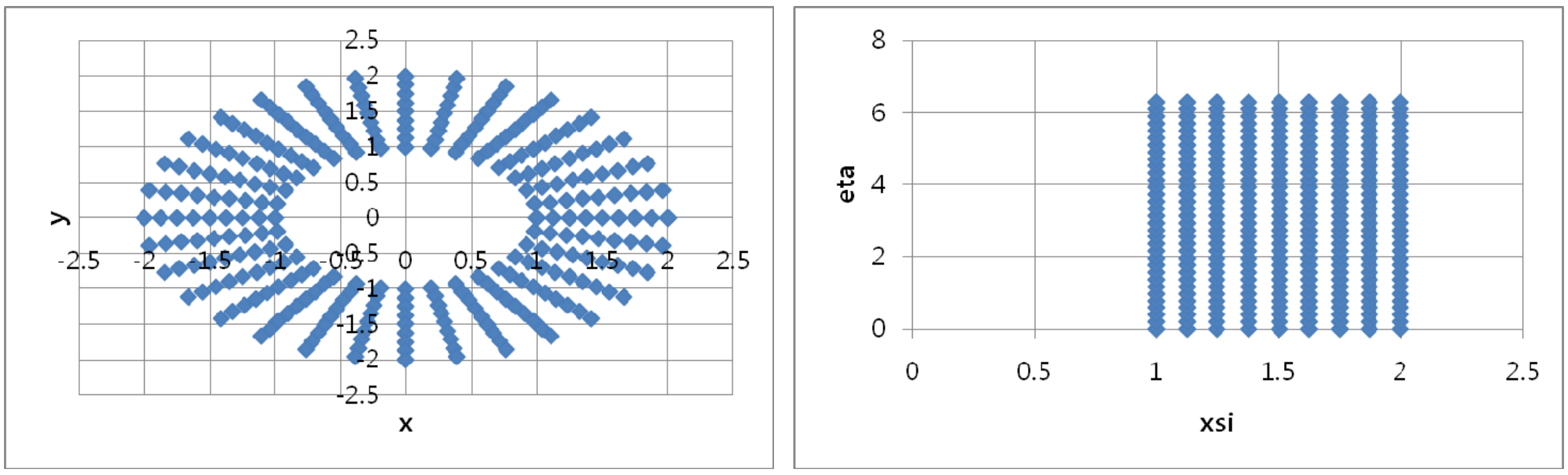

(a) Physical plane $(x,y)$ (b) Mapped plane $(\xi,\eta)$

Fig. 12. Generalization ability of normal mapping from physical plane to mapped plane ($I$=4, $J$=16, $M$=7).

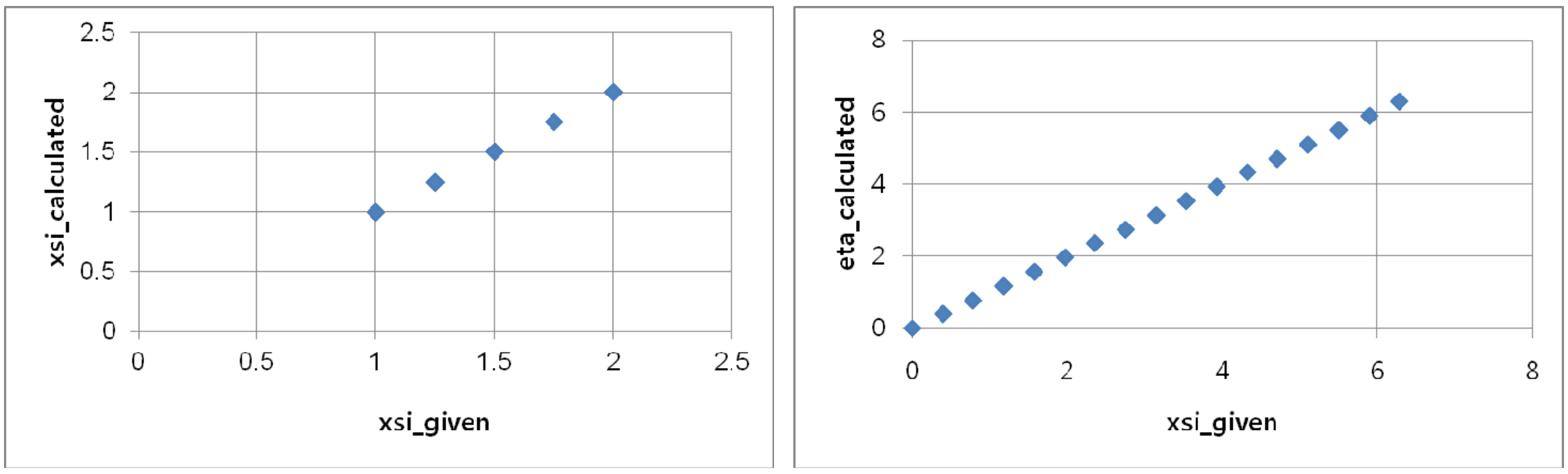

(a) $\xi$ (b) $\eta$

Fig. 13. Accuracy of normal mapping ($I$=4, $J$=16, $M$=7).

The solution of the inverse mapping $A'_{mn}$ and $B'_{mn}$ are determined by Eqs. (7) and (8). The mapping results are shown in Figs. 14 and 15. The accuracy of the inverse mapping is high as shown in Fig. 16. However, the generalization ability is not good as shown in Fig. 15.

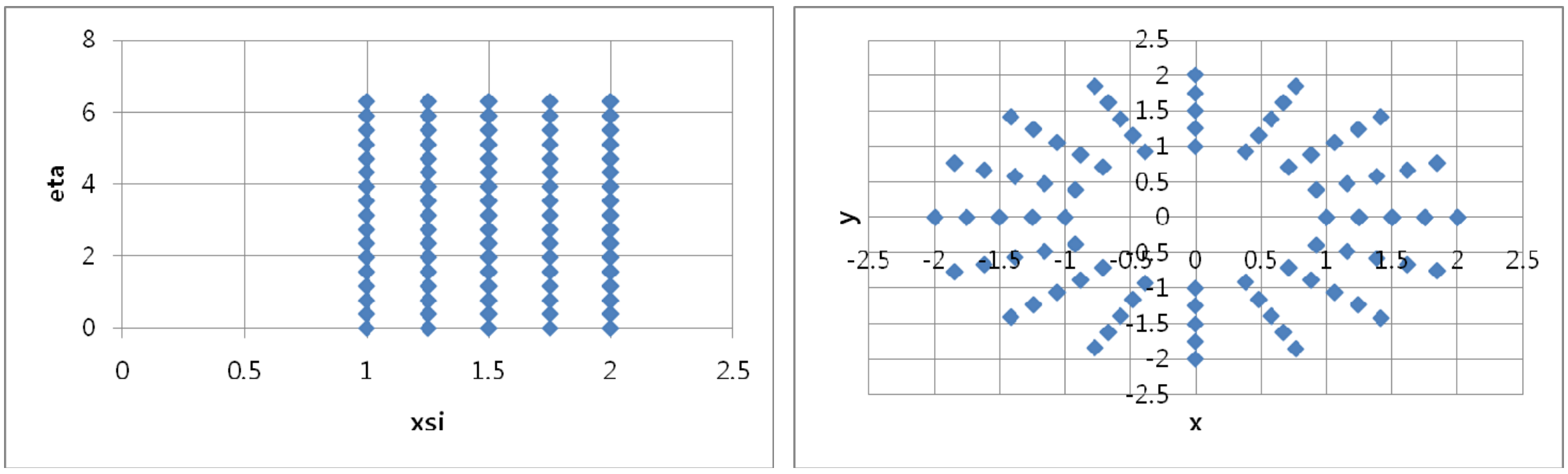

(a) Mapped plane $(\xi,\eta)$ (b) Physical plane $(x,y)$

Fig. 14. Inverse mapping $x=x(\xi,\eta)$, $y=y(\xi,\eta)$ ($I$=4, $J$=16, $M$=7).

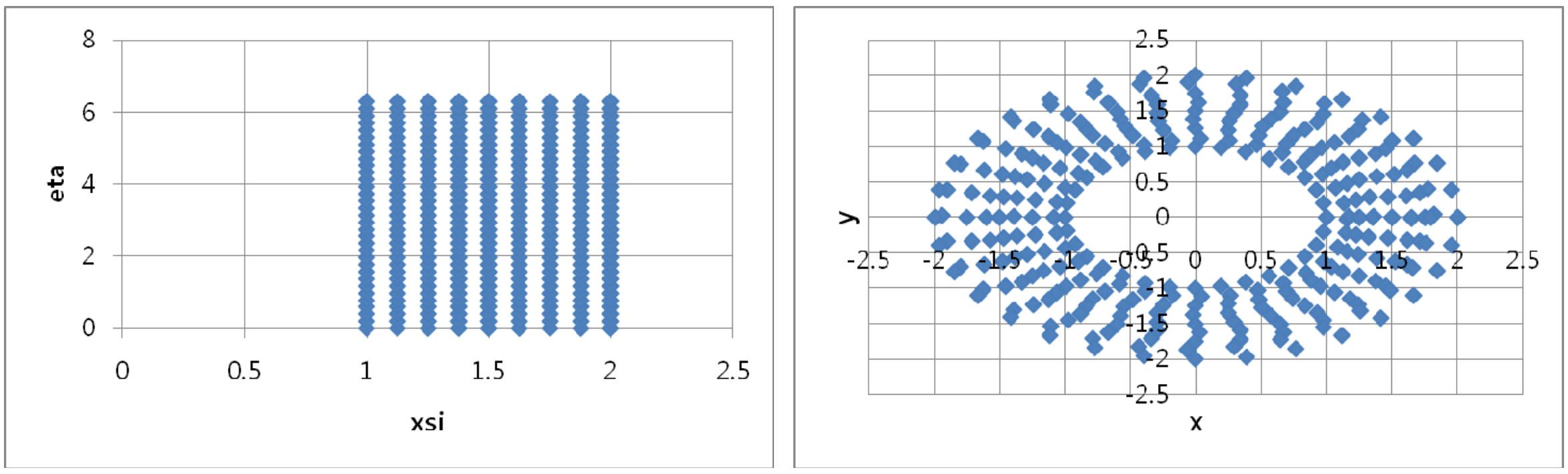

(a) Mapped plane $(\xi,\eta)$ (b) Physical plane $(x,y)$

Fig. 15. Generalization of physical plane $x=x(\xi,\eta)$, $y=y(\xi,\eta)$ ($I$=4, $J$=16, $M$=7).

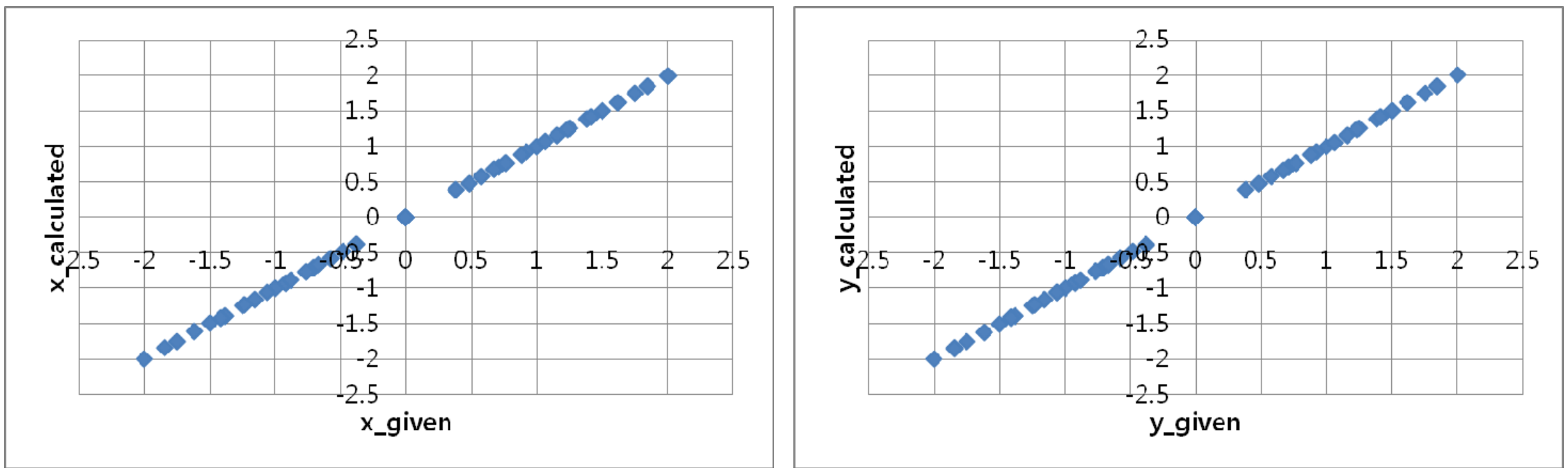

(a) $x$ (b) $y$

Fig. 16. Accuracy of inverse mappng ($M$=3) ($I$=4, $J$=16, $M$=7).

If we use the following parameters:

$$I = 5\,,\quad J = 20\,,\quad M = 7\,, \tag{13}$$

very good results are obtained not only for normal but also for inverse mappings. The results of the normal mapping are shown in Figs. (17) and (18). The accuracy and the generalization ability seem satisfactory.

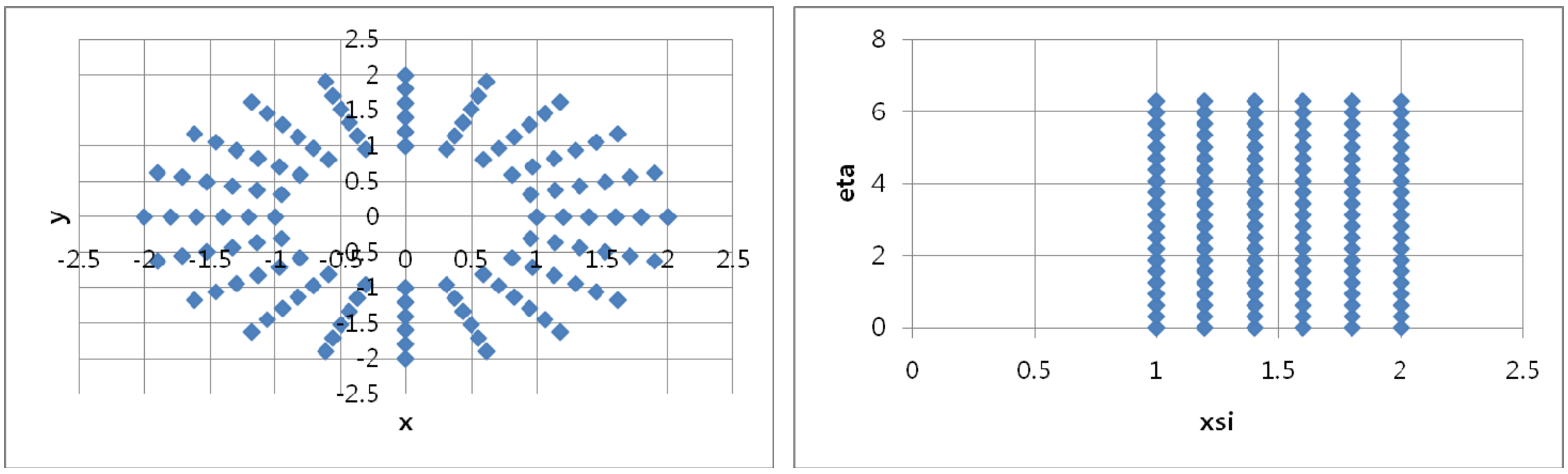

(a) Physical plane ($x$,$y$) (b) Mapped plane ($\xi$,$\eta$)

Fig. 17. Normal mapping $\xi=\xi(x,y)$, $\eta=\eta(x,y)$ ($I$=5, $J$=20, $M$=7).

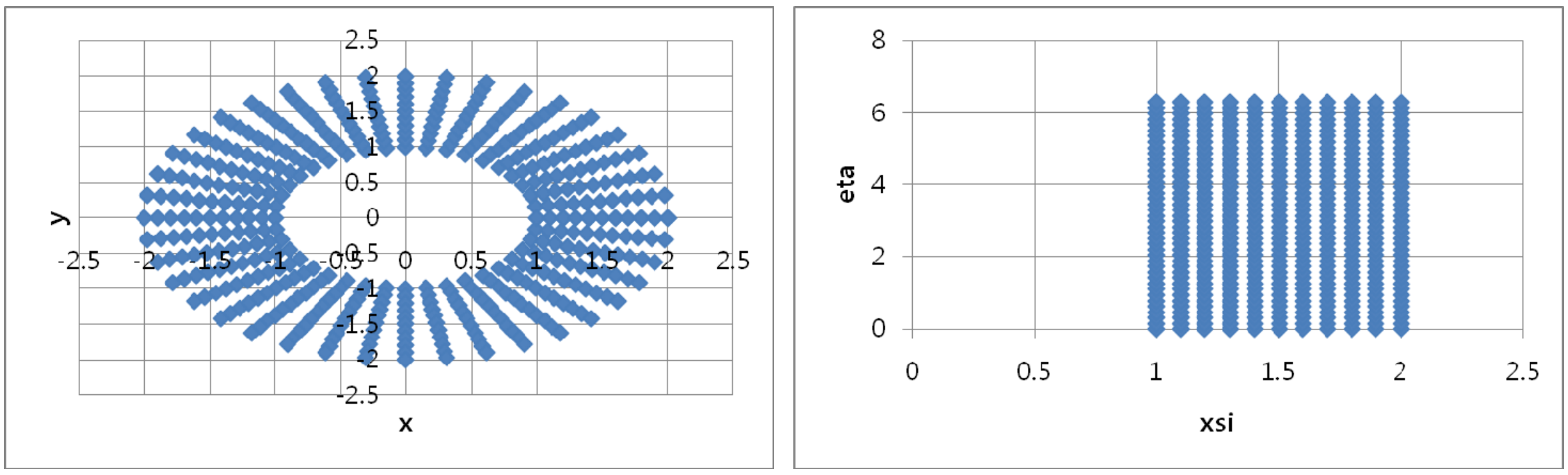

Fig. 18. Generalization ability of normal mapping from physical plane to mapped plane ($I$=5, $J$=20, $M$=7).

The results of the inverse mapping are shown in Figs. (19) and (20). The accuracy and the generalization ability seem satisfactory.

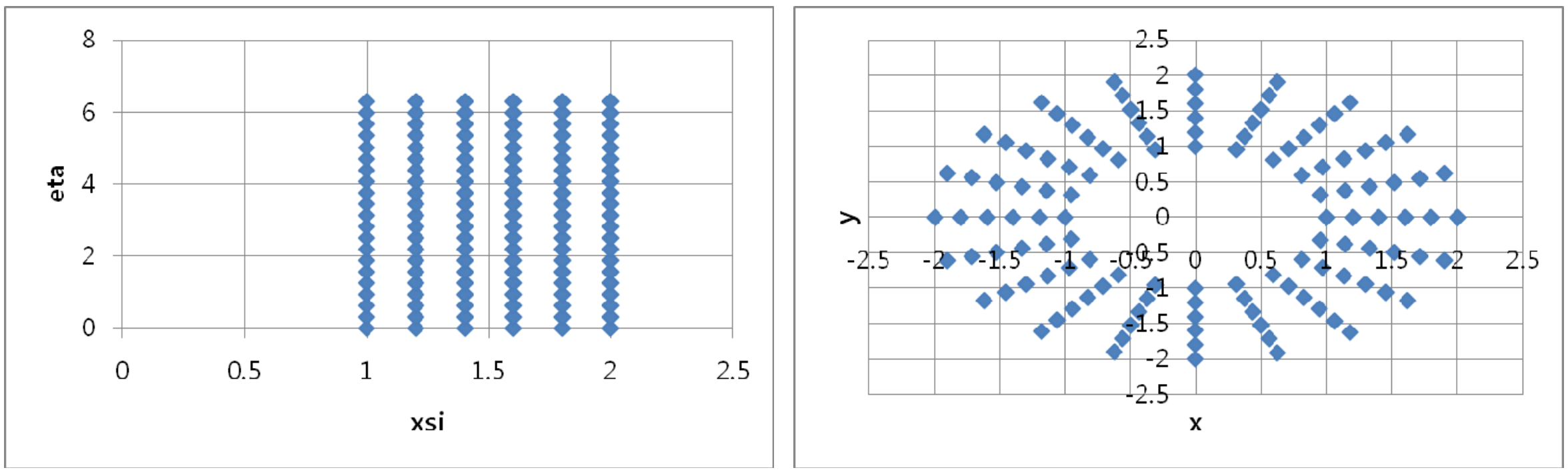

(a) Mapped plane $(\xi,\eta)$ (b) Physical plane $(x,y)$

Fig. 19. Inverse mapping $x=x(\xi,\eta)$, $y=y(\xi,\eta)$ ($I$=5, $J$=20, $M$=7).

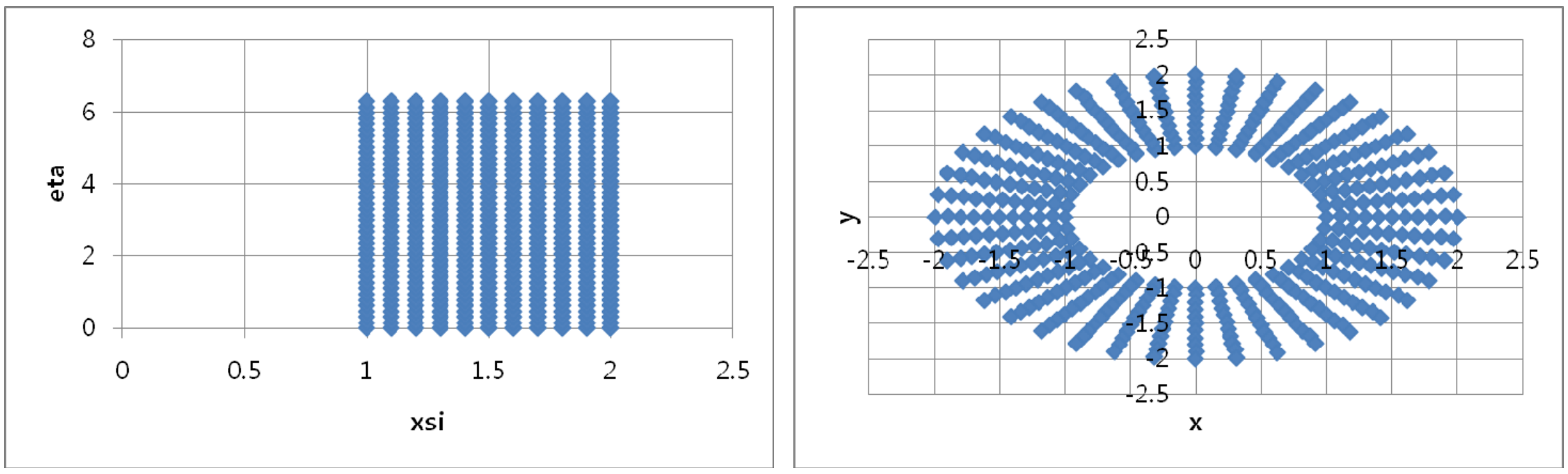

(a) Mapped plane $(\xi,\eta)$ (b) Physical plane $(x,y)$

Fig. 20. Generalization ability of inverse mapping from mapped plane to physical plane ($I$=5, $J$=20, $M$=7).

(2) Mapping of a 270deg circular sector

Given points for the rectangular coordinates in the physical plane are specified as

$$x_{ij} = \xi_i \cos\eta_j, \quad y_{ij} = \xi_i \sin\eta_j, \quad i = 0,1,\cdots,I, \quad j = 0,1,\cdots,J. \tag{14a}$$

And given points for the polar coordinates in the mapped plane are specified as

$$\xi_i = 1 + \frac{2-1}{I} i, \quad i = 0,1,\cdots,I; \qquad \eta_j = 0 + \frac{\pi - 0}{J} j, \quad j = 0,1,\cdots,J. \tag{14b}$$

In the calculations below, the following parameters are used:

$$I = 5, \quad J = 15, \quad M = 5. \tag{15}$$

The solution of the normal mapping $A_{mn}$ and $B_{mn}$ are determined by Eqs. (4) and (5). The mapping results are shown in Figs. 21 and 22.

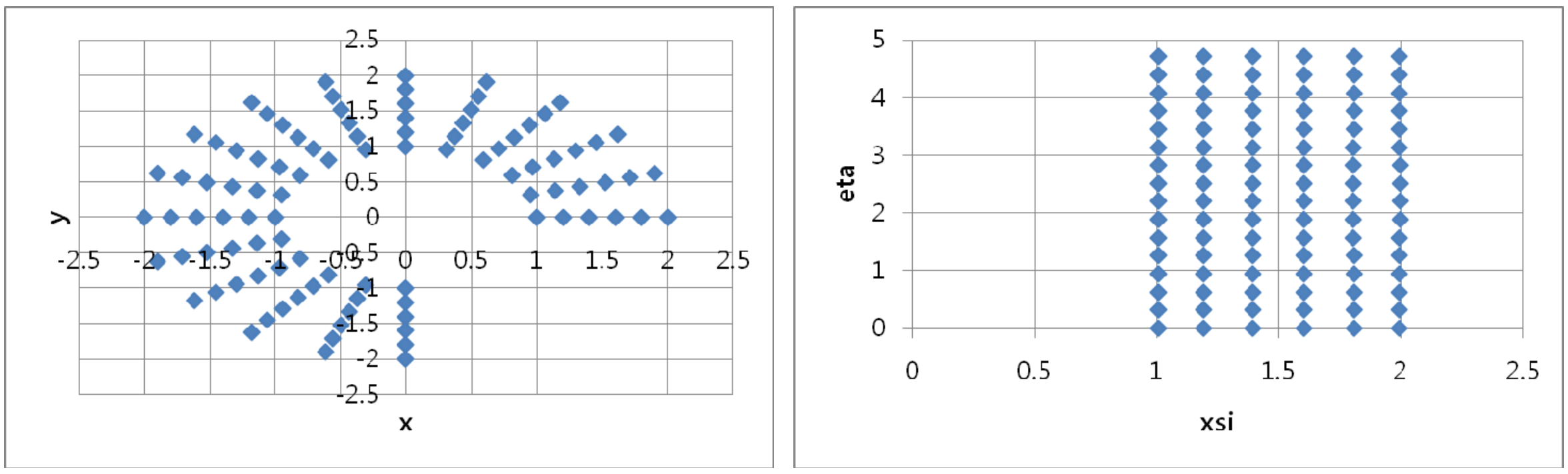

(a) Physical plane $(x,y)$ (b) Mapped plane $(\xi,\eta)$

Fig. 21. Normal mapping $\xi = \xi(x,y)$, $\eta = \eta(x,y)$ ($I$=5, $J$=15, $M$=5).

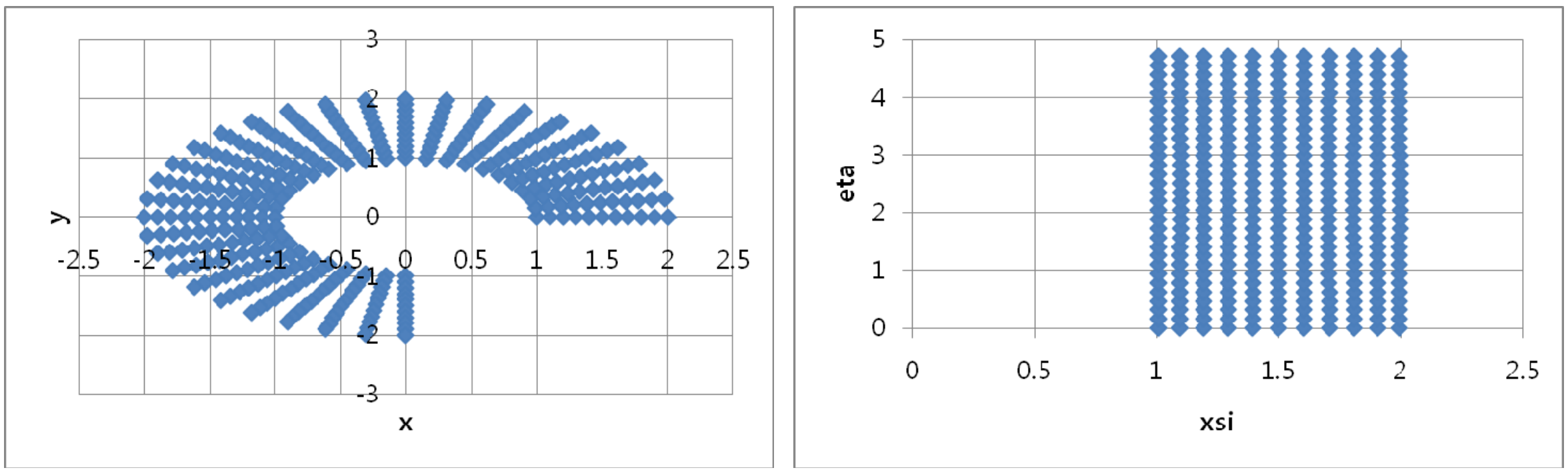

(a) Physical plane $(x,y)$ (b) Mapped plane $(\xi,\eta)$

Fig. 22. Generalization ability of normal mapping from physical plane to mapped plane ($I$=5, $J$=15, $M$=5).

The solution of the inverse mapping $A'_{mn}$ and $B'_{mn}$ are determined by Eqs. (7) and (8). The mapping results are shown in Figs. 23 and 24.

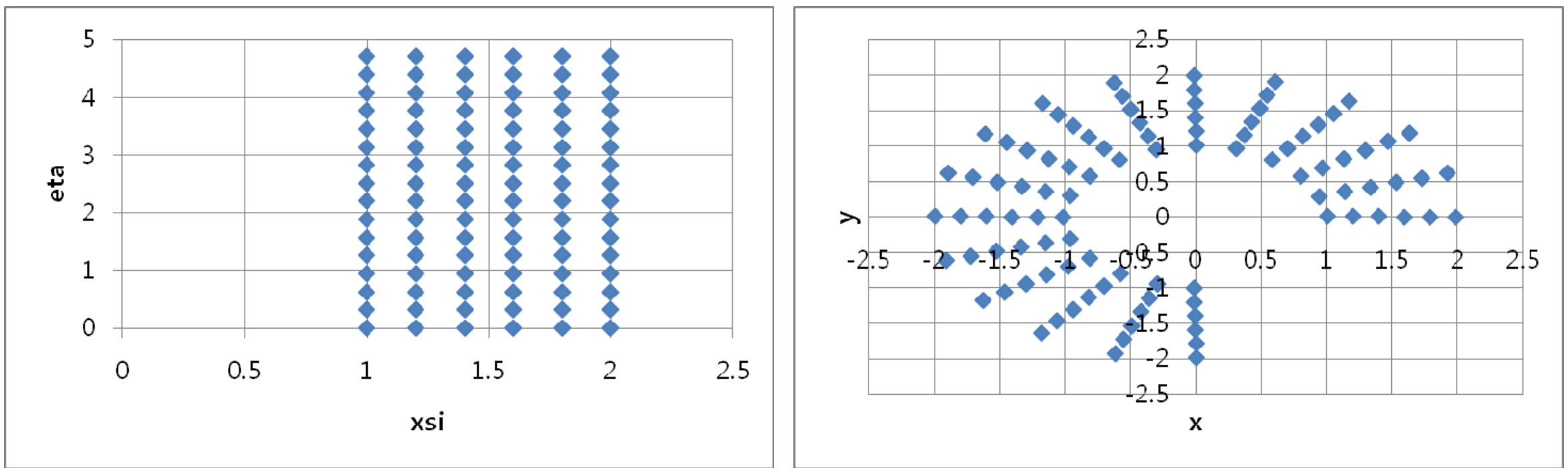

(a) Mapped plane $(\xi,\eta)$ (b) Physical plane $(x,y)$

Fig. 23. Inverse mapping $x=x(\xi,\eta)$, $y=y(\xi,\eta)$ ($I$=5, $J$=15, $M$=5).

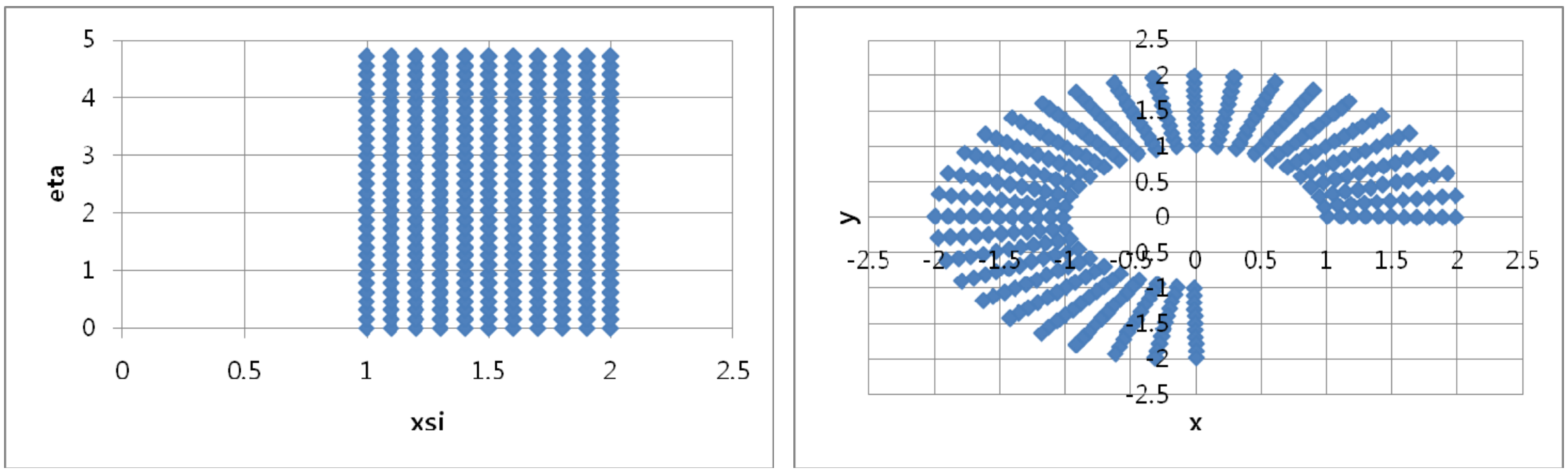

(a) Mapped plane $(\xi,\eta)$ (b) Physical plane $(x,y)$

Fig. 24. Generalization ability of inverse mapping from mapped plane to physical plane ($I$=5, $J$=15, $M$=5).

## 5. Application to Dirichlet problem of potential

We consider an application to a boundary value problem in an eccentric region in 2D as shown in Fig. 50. The inner boundary is a circle of radius $a$ with center at $(0,0)$. The outer boundary is a circle of radius $R$ with center at $(0, c_I)$. We define the $i$-th inner circle of radius $\xi_i$ with center at $(0, c_i)$, where

$$c_i = i\frac{c_I}{I}, \quad i = 0, 1, \cdots, I. \tag{16}$$

The $i$-th inner circle intersects the $x$-axis at $e_i^{\pm} = \xi_i \pm c_i$. We then have

$$e_i^- = -a - i\frac{R - c_I - a}{I}, \quad e_i^+ = a + i\frac{R + c_I - a}{I}, \tag{17}$$

$$c_i = \frac{1}{2}(e_i^+ + e_i^-) = i\frac{c_I}{I}, \quad \xi_i = \frac{1}{2}(e_i^+ - e_i^-) = a + i\frac{R-a}{I}. \tag{18}$$

We introduce the angular coordinates $\eta_j$:

$$\eta_j = \pi/J, \quad j = 0, 1, \cdots, J \tag{19}$$

along the circumferential direction. The Cartesian coordinates $(x_{ij}, y_{ij})$ corresponding to the curvilinear coordinates $(\xi_i, \eta_j)$ is given by

$$x_{ij} = c_i + \xi_i \cos\eta_j, \quad y_{ij} = \xi_i \sin\eta_j, \tag{20}$$

or

$$\xi_i = \sqrt{(x_{ij} - c_i)^2 + y_{ij}{}^2}, \quad \eta_j = \tan^{-1}\frac{y_{ij}}{x_{ij} - c_i}. \tag{21}$$

Using Eq. (18), Eq. (20) can be modified as

$$x_{ij} = c_i + \xi_i \cos\eta_j = \frac{c_I}{R-a}(\xi_i - a) + \xi_i \cos\eta_j. \tag{22}$$

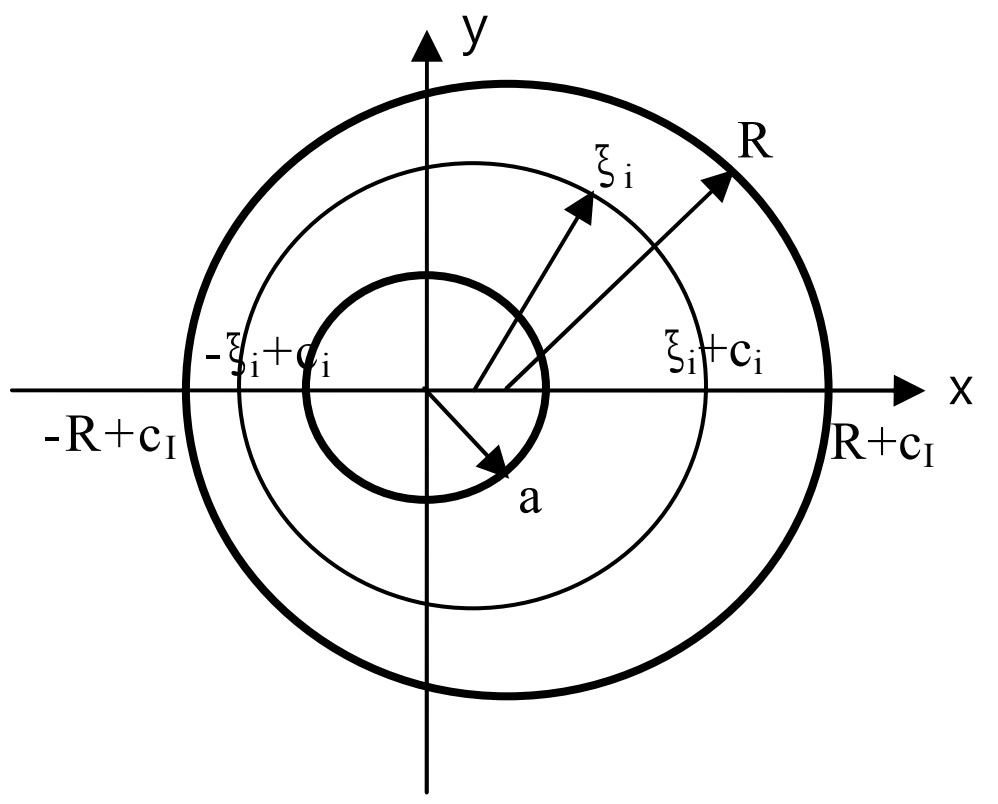

Fig. 25. Region surrounded by two eccentric circles with radius *R* and *a*.

The boundary value problem is a Dirichlet problem of potential and is defined in Cartesian coordinates as

$$\frac{\partial^2\phi}{\partial x^2}+\frac{\partial^2\phi}{\partial y^2}=0 \quad \text{in region surrounded by two eccentric circles,} \tag{23}$$

$$\phi=\phi_a \quad \text{on inner circle of radius } a\,, \tag{24}$$

$$\phi=\phi_R \quad \text{on inner circle of radius } R\,, \tag{25}$$

$$\frac{\partial\phi}{\partial y}=0 \quad \text{on } -R+c_I<x<R+c_I\,,\ \ y=0\,. \tag{26}$$

In Cartesian coordinates, the boundary value problem is written as

$$\left[\left(\frac{\partial\xi}{\partial x}\right)^2+\left(\frac{\partial\xi}{\partial y}\right)^2\right]\frac{\partial^2\phi}{\partial\xi^2}+2\left(\frac{\partial\xi}{\partial x}\frac{\partial\eta}{\partial x}+\frac{\partial\xi}{\partial y}\frac{\partial\eta}{\partial y}\right)\frac{\partial^2\phi}{\partial\xi\partial\eta}+\left[\left(\frac{\partial\eta}{\partial x}\right)^2+\left(\frac{\partial\eta}{\partial y}\right)^2\right]\frac{\partial^2\phi}{\partial\eta^2}, \\ +\left(\frac{\partial^2\xi}{\partial x^2}+\frac{\partial^2\xi}{\partial y^2}\right)\frac{\partial\phi}{\partial\xi}+\left(\frac{\partial^2\eta}{\partial x^2}+\frac{\partial^2\eta}{\partial y^2}\right)\frac{\partial\phi}{\partial\eta}=0. \tag{27}$$

$$\phi=\phi_a \quad \text{on } \xi=a\,,\ \ 0<\eta<\pi\,, \tag{28}$$

$$\phi=\phi_R \quad \text{on } \xi=R\,,\ \ 0<\eta<\pi\,, \tag{29}$$

$$\frac{\partial\phi}{\partial\eta}=0 \quad \text{on } a<\xi<R\,,\ \ \eta=0\,. \tag{30}$$

$$\frac{\partial\phi}{\partial\eta}=0 \quad \text{on } a<\xi<R\,,\ \ \eta=\pi\,. \tag{31}$$

The formulas of coordinates transformation are given in Eqs. (53)~(70) in Ref. [1]. The numerical inverse mapping $x=x(\xi,\eta)$ and $y=y(\xi,\eta)$ can be conducted using the inverse transformation given by Eqs. (7)~(8).

(1) Region surrounded by two concentric circles

The first example is a Dirichlet problem in a region surrounded by two concentric circles. We used the following parameters for the numerical calculation:

$$a=2\,,\ \ R=6\,,\ \ c_I=0\,,\ \ I=4\,,\ \ J=6\,,\ \ M=6\,,\ \ \phi_a=0\,,\ \ \phi_R=1\,. \tag{32}$$

The coefficients in the inverse mapping given by Eq. (6) is given in Table 4. The numerical The numerical inverse mapping $x=x(\xi,\eta)$ and $y=y(\xi,\eta)$ is shown in Fig. 26.

Table 4. Coeffficients of interpolation functions in Eq. (6).

| m | n | A’ | B’ | m | n | A’ | B’ |
|---|---|---|---|---|---|---|---|
| 0 | 0 | 0.130265 | 4.008964 | 5 | 0 | 0.001289 | 0.011132 |
| 1 | 0 | 0.864018 | -5.2119 | 5 | 1 | 0.000421 | -0.02484 |
| 1 | 1 | 0.015173 | -0.83887 | 5 | 2 | 0 | 0 |
| 2 | 0 | 0.028228 | 2.357771 | 5 | 3 | 0 | 0 |
| 2 | 1 | -0.0247 | 2.283302 | 5 | 4 | 0.059702 | 0.036893 |
| 2 | 2 | 0 | -0.21154 | 5 | 5 | 0 | 0.048884 |
| 3 | 0 | 0.014572 | -0.37876 | 6 | 0 | -7.3E-05 | -0.00084 |
| 3 | 1 | 0.012222 | -0.72033 | 6 | 1 | -2.1E-05 | 0.001242 |
| 3 | 2 | -0.48751 | 0.05004 | 6 | 2 | 0 | 0 |
| 3 | 3 | 0 | 0.284279 | 6 | 3 | 0 | 0 |
| 4 | 0 | -0.00777 | -0.01806 | 6 | 4 | 0 | 0 |
| 4 | 1 | -0.00327 | 0.192502 | 6 | 5 | -0.0076 | 0 |
| 4 | 2 | 0 | 0 | 6 | 6 | 0 | -0.00519 |
| 4 | 3 | -0.02159 | -0.23181 | | | | |
| 4 | 4 | 0 | -0.17322 | | | | |

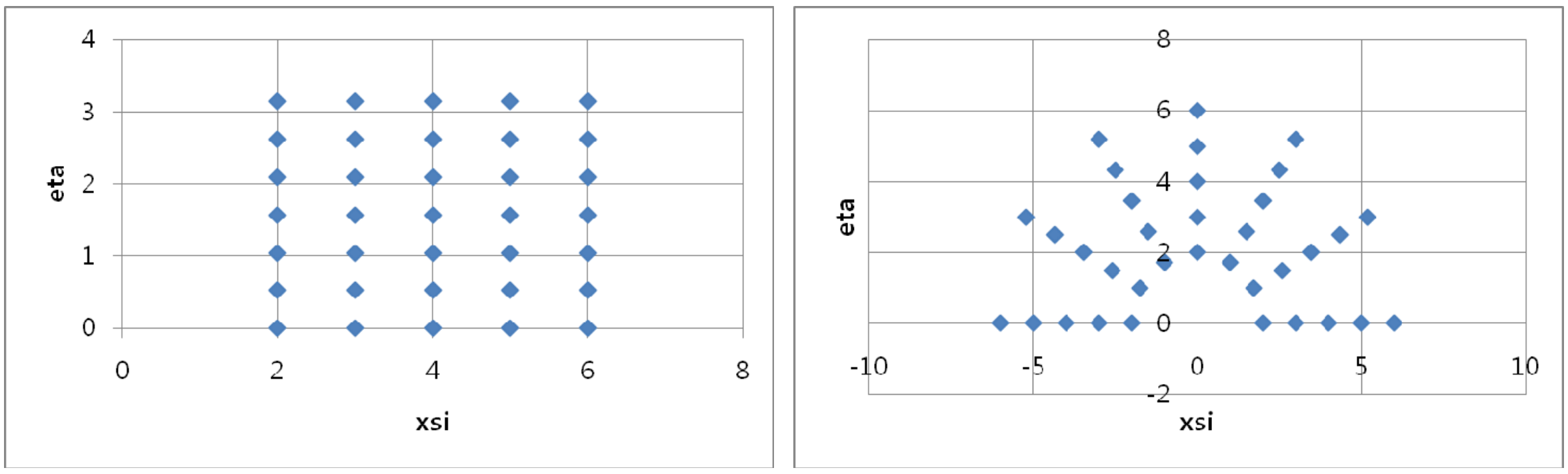

(a) Mapped plane $(\xi,\eta)$ (b) Physical plane $(x,y)$

Fig. 26. Mapping of region surrounded by two concentric circles.

The numerical and exact solution of $\phi$ is shown in Table 10. The agreement between the numerical and exact solutions is satisfactory, where the exact solution is given by

$$\phi = \frac{\phi_R \ln(r/a) - \phi_a \ln(r/R)}{\ln(R/a)}. \tag{33}$$

Table 5. Numerical and exact solutions of potential $\phi$.

| $\xi$ | $\eta$ = 0deg | $\eta$ = 30 | $\eta$ = 60 | $\eta$ = 90 | $\eta$ =120 | $\eta$ =150 | $\eta$ =180 | exact |
|---|---|---|---|---|---|---|---|---|
| 2 | 0 | 0 | 0 | 0 | 0 | 0 | 0 | 0 |
| 3 | 0.375004 | 0.375004 | 0.376727 | 0.376666 | 0.374204 | 0.371361 | 0.371361 | 0.36907 |
| 4 | 0.630386 | 0.630386 | 0.629934 | 0.629739 | 0.629603 | 0.62986 | 0.62986 | 0.63093 |
| 5 | 0.829118 | 0.829118 | 0.826957 | 0.826588 | 0.82811 | 0.830754 | 0.830754 | 0.834044 |
| 6 | 1 | 1 | 1 | 1 | 1 | 1 | 1 | 1 |

(2) Region surrounded by two eccentric circles

The second example is a Dirichlet problem in a region surrounded by two eccentric circles. We used the following parameters for the numerical calculation:

$$a = 2,\quad R = 6,\quad c_I = 2,\quad I = 4,\quad J = 6,\quad M = 6,\quad \phi_a = 0,\quad \phi_R = 1. \tag{34}$$

The coefficients in the inverse mapping given by Eq. (6) is given in Table 12. The numerical inverse mapping $x = x(\xi,\eta)$ and $y = y(\xi,\eta)$ is shown in Fig. 27.

Table 6. Coeffficients of interpolation functions in Eq. (6).

| m | n | Ap | Bp | m | n | Ap | Bp |
|---|---|---|---|---|---|---|---|
| 0 | 0 | -1.04155 | 4.008964 | 5 | 0 | 0.001477 | 0.011132 |
| 1 | 0 | 1.611337 | -5.2119 | 5 | 1 | -0.00217 | -0.02484 |
| 1 | 1 | -0.07819 | -0.83887 | 5 | 2 | 0 | 0 |
| 2 | 0 | -0.10756 | 2.357771 | 5 | 3 | 0 | 0 |
| 2 | 1 | 0.110674 | 2.283302 | 5 | 4 | 0.059702 | 0.036893 |
| 2 | 2 | 0 | -0.21154 | 5 | 5 | 0 | 0.048884 |
| 3 | 0 | 0.050103 | -0.37876 | 6 | 0 | -7.1E-05 | -0.00084 |
| 3 | 1 | -0.06298 | -0.72033 | 6 | 1 | 0.000109 | 0.001242 |
| 3 | 2 | -0.48751 | 0.05004 | 6 | 2 | 0 | 0 |
| 3 | 3 | 0 | 0.284279 | 6 | 3 | 0 | 0 |
| 4 | 0 | -0.01215 | -0.01806 | 6 | 4 | 0 | 0 |
| 4 | 1 | 0.016832 | 0.192502 | 6 | 5 | -0.0076 | 0 |
| 4 | 2 | 0 | 0 | 6 | 6 | 0 | -0.00519 |
| 4 | 3 | -0.02159 | -0.23181 | | | | |
| 4 | 4 | 0 | -0.17322 | | | | |

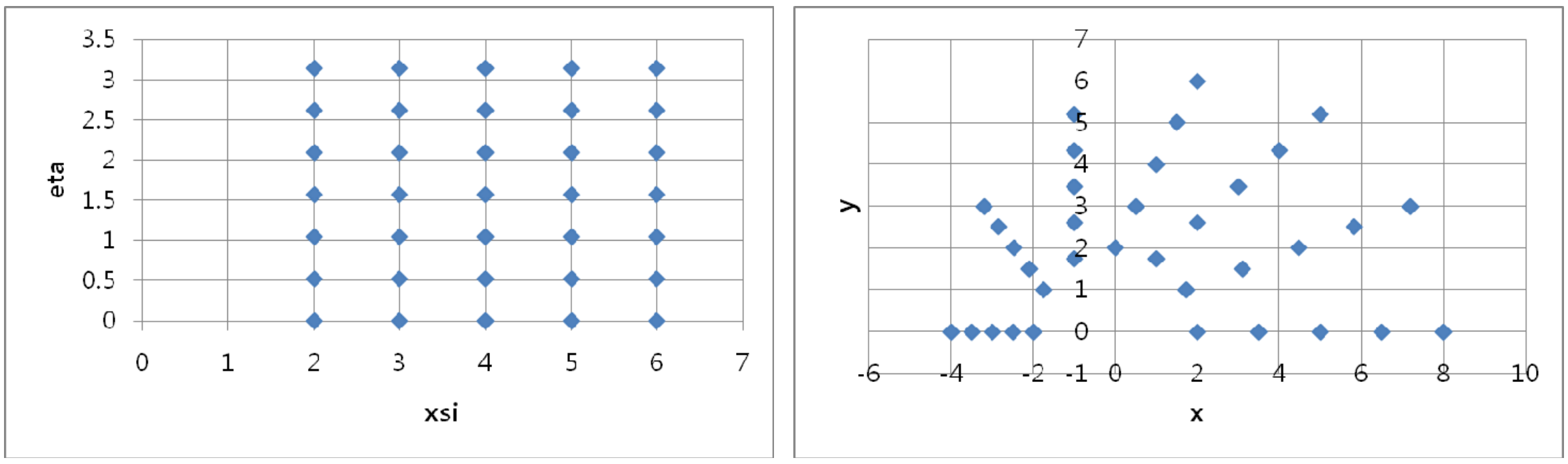

(a) Mapped plane $(\xi,\eta)$ (b) Physical plane $(x,y)$

Fig. 27. Mapping of region surrounded by two eccentric circles.

The accuracy of mapping is shown in Tables 7~18 and seems satisfactory.

Table 7. Accuracy of $x$=$x(\xi,\eta)$.

| $\eta$ | $\xi$= 2.0 | | $\xi$= 3.0 | | $\xi$= 4.0 | | $\xi$= 5.0 | | $\xi$= 6.0 | |
|---|---|---|---|---|---|---|---|---|---|---|
| deg | *Num.* | *Exact* | *Num.* | *Exact* | *Num.* | *Exact* | *Num.* | *Exact* | *Num.* | *Exact* |
| 0 | 2 | 2 | 3.5 | 3.5 | 5 | 5 | 6.5 | 6.5 | 8 | 8 |
| 30 | 1.732051 | 1.732051 | 3.098076 | 3.098076 | 4.464102 | 4.464102 | 5.830127 | 5.830127 | 7.196152 | 7.196152 |
| 60 | 1 | 1 | 2 | 2 | 3 | 3 | 4 | 4 | 5 | 5 |
| 90 | 0 | 0 | 0.5 | 0.5 | 1 | 1 | 1.5 | 1.5 | 2 | 2 |
| 120 | -1 | -1 | -1 | -1 | -1 | -1 | -1 | -1 | -1 | -1 |
| 150 | -1.73205 | -1.73205 | -2.09808 | -2.09808 | -2.4641 | -2.4641 | -2.83013 | -2.83013 | -3.19615 | -3.19615 |
| 180 | -2 | -2 | -2.5 | -2.5 | -3 | -3 | -3.5 | -3.5 | -4 | -4 |

Table 8. Accuracy of $x_\xi$=$x_\xi(\xi,\eta)$.

| $\eta$ | $\xi$= 2.0 | | $\xi$= 3.0 | | $\xi$= 4.0 | | $\xi$= 5.0 | | $\xi$= 6.0 | |
|---|---|---|---|---|---|---|---|---|---|---|
| deg | *Num.* | *Exact* | *Num.* | *Exact* | *Num.* | *Exact* | *Num.* | *Exact* | *Num.* | *Exact* |
| 0 | 1.497979 | 1.5 | 1.500931 | 1.5 | 1.499095 | 1.5 | 1.501783 | 1.5 | 1.491166 | 1.5 |
| 30 | 1.365369 | 1.366025 | 1.366615 | 1.366025 | 1.365348 | 1.366025 | 1.367467 | 1.366025 | 1.358556 | 1.366025 |
| 60 | 1.000708 | 1 | 1.000249 | 1 | 0.99955 | 1 | 1.0011 | 1 | 0.993895 | 1 |
| 90 | 0.502072 | 0.5 | 0.499908 | 0.5 | 0.499778 | 0.5 | 0.500759 | 0.5 | 0.49526 | 0.5 |
| 120 | 0.003437 | 0 | -0.00043 | 0 | 0.000005 | 0 | 0.000418 | 0 | -0.00338 | 0 |
| 150 | -0.36122 | -0.36603 | -0.3668 | -0.36603 | -0.36579 | -0.36603 | -0.36595 | -0.36603 | -0.36804 | -0.36603 |
| 180 | -0.49383 | -0.5 | -0.50112 | -0.5 | -0.49954 | -0.5 | -0.50026 | -0.5 | -0.50065 | -0.5 |

Table 9. Accuracy of $x_\eta = x_\eta(\xi, \eta)$.

| $\eta$ deg | $\xi$= 2.0 Num. | $\xi$= 2.0 Exact | $\xi$= 3.0 Num. | $\xi$= 3.0 Exact | $\xi$= 4.0 Num. | $\xi$= 4.0 Exact | $\xi$= 5.0 Num. | $\xi$= 5.0 Exact | $\xi$= 6.0 Num. | $\xi$= 6.0 Exact |
|---|---|---|---|---|---|---|---|---|---|---|
| 0 | -0.00539 | 0 | -0.00809 | 0 | -0.01078 | 0 | -0.01348 | 0 | -0.01617 | 0 |
| 30 | -0.99908 | -1 | -1.49863 | -1.5 | -1.99817 | -2 | -2.49771 | -2.5 | -2.99725 | -3 |
| 60 | -1.73242 | -1.73205 | -2.59863 | -2.59808 | -3.46484 | -3.4641 | -4.33105 | -4.33013 | -5.19726 | -5.19615 |
| 90 | -1.99972 | -2 | -2.99958 | -3 | -3.99944 | -4 | -4.9993 | -5 | -5.99916 | -6 |
| 120 | -1.73242 | -1.73205 | -2.59863 | -2.59808 | -3.46484 | -3.4641 | -4.33105 | -4.33013 | -5.19726 | -5.19615 |
| 150 | -0.99908 | -1 | -1.49863 | -1.5 | -1.99817 | -2 | -2.49771 | -2.5 | -2.99725 | -3 |
| 180 | -0.00539 | 0 | -0.00809 | 0 | -0.01078 | 0 | -0.01348 | 0 | -0.01617 | 0 |

Table 10. Accuracy of $x_{\xi\xi} = x_{\xi\xi}(\xi, \eta)$.

| $\eta$ deg | $\xi$= 2.0 Num. | $\xi$= 2.0 Exact | $\xi$= 3.0 Num. | $\xi$= 3.0 Exact | $\xi$= 4.0 Num. | $\xi$= 4.0 Exact | $\xi$= 5.0 Num. | $\xi$= 5.0 Exact | $\xi$= 6.0 Num. | $\xi$= 6.0 Exact |
|---|---|---|---|---|---|---|---|---|---|---|
| 0 | 0.005016 | 0 | -0.0007 | 0 | -0.00057 | 0 | 0.003823 | 0 | -0.04021 | 0 |
| 30 | -0.00067 | 0 | -0.00013 | 0 | -0.00057 | 0 | 0.003254 | 0 | -0.03453 | 0 |
| 60 | -0.00636 | 0 | 0.000437 | 0 | -0.00057 | 0 | 0.002686 | 0 | -0.02884 | 0 |
| 90 | -0.01204 | 0 | 0.001005 | 0 | -0.00057 | 0 | 0.002117 | 0 | -0.02316 | 0 |
| 120 | -0.01773 | 0 | 0.001574 | 0 | -0.00057 | 0 | 0.001548 | 0 | -0.01747 | 0 |
| 150 | -0.02341 | 0 | 0.002143 | 0 | -0.00057 | 0 | 0.00098 | 0 | -0.01179 | 0 |
| 180 | -0.0291 | 0 | 0.002711 | 0 | -0.00057 | 0 | 0.000411 | 0 | -0.0061 | 0 |

Table 11. Accuracy of $x_{\xi\eta} = x_{\xi\eta}(\xi, \eta)$.

| $\eta$ deg | $\xi$= 2.0 Num. | $\xi$= 2.0 Exact | $\xi$= 3.0 Num. | $\xi$= 3.0 Exact | $\xi$= 4.0 Num. | $\xi$= 4.0 Exact | $\xi$= 5.0 Num. | $\xi$= 5.0 Exact | $\xi$= 6.0 Num. | $\xi$= 6.0 Exact |
|---|---|---|---|---|---|---|---|---|---|---|
| 0 | -8.9E-05 | 0 | -0.00335 | 0 | -0.00226 | 0 | -0.00335 | 0 | -8.9E-05 | 0 |
| 30 | -0.49694 | -0.5 | -0.50019 | -0.5 | -0.49911 | -0.5 | -0.50019 | -0.5 | -0.49694 | -0.5 |
| 60 | -0.86361 | -0.86603 | -0.86686 | -0.86603 | -0.86578 | -0.86603 | -0.86686 | -0.86603 | -0.86361 | -0.86603 |
| 90 | -0.99725 | -1 | -1.00051 | -1 | -0.99943 | -1 | -1.00051 | -1 | -0.99725 | -1 |
| 120 | -0.86361 | -0.86603 | -0.86686 | -0.86603 | -0.86578 | -0.86603 | -0.86686 | -0.86603 | -0.86361 | -0.86603 |
| 150 | -0.49694 | -0.5 | -0.50019 | -0.5 | -0.49911 | -0.5 | -0.50019 | -0.5 | -0.49694 | -0.5 |
| 180 | -8.9E-05 | 0 | -0.00335 | 0 | -0.00226 | 0 | -0.00335 | 0 | -8.9E-05 | 0 |

Table 12. Accuracy of $x_{\eta\eta} = x_{\eta\eta}(\xi, \eta)$.

| $\eta$ deg | $\xi$= 2.0 Num. | $\xi$= 2.0 Exact | $\xi$= 3.0 Num. | $\xi$= 3.0 Exact | $\xi$= 4.0 Num. | $\xi$= 4.0 Exact | $\xi$= 5.0 Num. | $\xi$= 5.0 Exact | $\xi$= 6.0 Num. | $\xi$= 6.0 Exact |
|---|---|---|---|---|---|---|---|---|---|---|
| 0 | -1.95002 | -2 | -2.92503 | -3 | -3.90005 | -4 | -4.87506 | -5 | -5.85007 | -6 |
| 30 | -1.73649 | -1.73205 | -2.60473 | -2.59808 | -3.47297 | -3.4641 | -4.34122 | -4.33013 | -5.20946 | -5.19615 |
| 60 | -0.99919 | -1 | -1.49878 | -1.5 | -1.99837 | -2 | -2.49796 | -2.5 | -2.99756 | -3 |
| 90 | 0 | 0 | 0 | 0 | 0 | 0 | 0 | 0 | 0 | 0 |
| 120 | 0.999185 | 1 | 1.498778 | 1.5 | 1.99837 | 2 | 2.497963 | 2.5 | 2.997555 | 3 |
| 150 | 1.736487 | 1.732051 | 2.60473 | 2.598076 | 3.472974 | 3.464102 | 4.341217 | 4.330127 | 5.209461 | 5.196152 |
| 180 | 1.950022 | 2 | 2.925033 | 3 | 3.900045 | 4 | 4.875056 | 5 | 5.850067 | 6 |

Table 13. Accuracy of $y = y(\xi, \eta)$.

| $\eta$ deg | $\xi$= 2.0 Num. | $\xi$= 2.0 Exact | $\xi$= 3.0 Num. | $\xi$= 3.0 Exact | $\xi$= 4.0 Num. | $\xi$= 4.0 Exact | $\xi$= 5.0 Num. | $\xi$= 5.0 Exact | $\xi$= 6.0 Num. | $\xi$= 6.0 Exact |
|---|---|---|---|---|---|---|---|---|---|---|
| 0 | -4.2E-05 | 0 | -2.1E-05 | 0 | 0 | 0 | 0.000021 | 0 | 0.000042 | 0 |
| 30 | 1.00025 | 1 | 1.500125 | 1.5 | 2 | 2 | 2.499875 | 2.5 | 2.99975 | 3 |
| 60 | 1.731426 | 1.732051 | 2.597764 | 2.598076 | 3.464102 | 3.464102 | 4.330439 | 4.330127 | 5.196777 | 5.196152 |
| 90 | 2.000833 | 2 | 3.000416 | 3 | 4 | 4 | 4.999584 | 5 | 5.999167 | 6 |
| 120 | 1.731426 | 1.732051 | 2.597764 | 2.598076 | 3.464102 | 3.464102 | 4.330439 | 4.330127 | 5.196777 | 5.196152 |
| 150 | 1.00025 | 1 | 1.500125 | 1.5 | 2 | 2 | 2.499875 | 2.5 | 2.99975 | 3 |
| 180 | -4.2E-05 | 0 | -2.1E-05 | 0 | 0 | 0 | 0.000021 | 0 | 0.000042 | 0 |

Table 14. Accuracy of $y_\xi = y_\xi(\xi, \eta)$.

| $\eta$ deg | $\xi$= 2.0 Num. | $\xi$= 2.0 Exact | $\xi$= 3.0 Num. | $\xi$= 3.0 Exact | $\xi$= 4.0 Num. | $\xi$= 4.0 Exact | $\xi$= 5.0 Num. | $\xi$= 5.0 Exact | $\xi$= 6.0 Num. | $\xi$= 6.0 Exact |
|---|---|---|---|---|---|---|---|---|---|---|
| 0 | -0.17369 | 0 | 0.04846 | 0 | -0.03561 | 0 | 0.05848 | 0 | -0.25385 | 0 |
| 30 | 0.341767 | 0.5 | 0.544412 | 0.5 | 0.466844 | 0.5 | 0.554432 | 0.5 | 0.261607 | 0.5 |
| 60 | 0.723836 | 0.866025 | 0.906973 | 0.866025 | 0.835908 | 0.866025 | 0.916993 | 0.866025 | 0.643677 | 0.866025 |
| 90 | 0.872689 | 1 | 1.036317 | 1 | 0.971755 | 1 | 1.046337 | 1 | 0.792529 | 1 |
| 120 | 0.75505 | 0.866025 | 0.89917 | 0.866025 | 0.84111 | 0.866025 | 0.90919 | 0.866025 | 0.67489 | 0.866025 |
| 150 | 0.404194 | 0.5 | 0.528805 | 0.5 | 0.477248 | 0.5 | 0.538825 | 0.5 | 0.324034 | 0.5 |
| 180 | -0.08005 | 0 | 0.025049 | 0 | -0.02001 | 0 | 0.035069 | 0 | -0.16021 | 0 |

Table 15. Accuracy of $y_\eta = y_\eta(\xi,\eta)$.

| $\eta$ deg | $\xi$= 2.0 Num. | $\xi$= 2.0 Exact | $\xi$= 3.0 Num. | $\xi$= 3.0 Exact | $\xi$= 4.0 Num. | $\xi$= 4.0 Exact | $\xi$= 5.0 Num. | $\xi$= 5.0 Exact | $\xi$= 6.0 Num. | $\xi$= 6.0 Exact |
|---|---|---|---|---|---|---|---|---|---|---|
| 0 | 2.028745 | 2 | 3.01545 | 3 | 4.002155 | 4 | 4.98886 | 5 | 5.975565 | 6 |
| 30 | 1.724581 | 1.732051 | 2.59422 | 2.598076 | 3.463858 | 3.464102 | 4.333497 | 4.330127 | 5.203136 | 5.196152 |
| 60 | 1.003309 | 1 | 1.501679 | 1.5 | 2.000049 | 2 | 2.498419 | 2.5 | 2.996789 | 3 |
| 90 | 0 | 0 | 0 | 0 | 0 | 0 | 0 | 0 | 0 | 0 |
| 120 | -1.00331 | -1 | -1.50168 | -1.5 | -2.00005 | -2 | -2.49842 | -2.5 | -2.99679 | -3 |
| 150 | -1.72458 | -1.73205 | -2.59422 | -2.59808 | -3.46386 | -3.4641 | -4.3335 | -4.33013 | -5.20314 | -5.19615 |
| 180 | -2.02875 | -2 | -3.01545 | -3 | -4.00216 | -4 | -4.98886 | -5 | -5.97557 | -6 |

Table 16. Accuracy of $y_{\xi\xi} = y_{\xi\xi}(\xi,\eta)$.

| $\eta$ deg | $\xi$= 2.0 Num. | $\xi$= 2.0 Exact | $\xi$= 3.0 Num. | $\xi$= 3.0 Exact | $\xi$= 4.0 Num. | $\xi$= 4.0 Exact | $\xi$= 5.0 Num. | $\xi$= 5.0 Exact | $\xi$= 6.0 Num. | $\xi$= 6.0 Exact |
|---|---|---|---|---|---|---|---|---|---|---|
| 0 | 0.683732 | 0 | -0.07071 | 0 | -0.00668 | 0 | 0.107451 | 0 | -1.09789 | 0 |
| 30 | 0.618704 | 0 | -0.06421 | 0 | -0.00668 | 0 | 0.100948 | 0 | -1.03286 | 0 |
| 60 | 0.553675 | 0 | -0.05771 | 0 | -0.00668 | 0 | 0.094445 | 0 | -0.96784 | 0 |
| 90 | 0.488647 | 0 | -0.0512 | 0 | -0.00668 | 0 | 0.087943 | 0 | -0.90281 | 0 |
| 120 | 0.423618 | 0 | -0.0447 | 0 | -0.00668 | 0 | 0.08144 | 0 | -0.83778 | 0 |
| 150 | 0.35859 | 0 | -0.0382 | 0 | -0.00668 | 0 | 0.074937 | 0 | -0.77275 | 0 |
| 180 | 0.293561 | 0 | -0.03169 | 0 | -0.00668 | 0 | 0.068434 | 0 | -0.70772 | 0 |

Table 17. Accuracy of $y_{\xi\eta} = y_{\xi\eta}(\xi,\eta)$.

| $\eta$ deg | $\xi$= 2.0 Num. | $\xi$= 2.0 Exact | $\xi$= 3.0 Num. | $\xi$= 3.0 Exact | $\xi$= 4.0 Num. | $\xi$= 4.0 Exact | $\xi$= 5.0 Num. | $\xi$= 5.0 Exact | $\xi$= 6.0 Num. | $\xi$= 6.0 Exact |
|---|---|---|---|---|---|---|---|---|---|---|
| 0 | 1.016512 | 1 | 0.979253 | 1 | 0.991673 | 1 | 0.979253 | 1 | 1.016512 | 1 |
| 30 | 0.899446 | 0.866025 | 0.862187 | 0.866025 | 0.874607 | 0.866025 | 0.862187 | 0.866025 | 0.899446 | 0.866025 |
| 60 | 0.528177 | 0.5 | 0.490918 | 0.5 | 0.503338 | 0.5 | 0.490918 | 0.5 | 0.528177 | 0.5 |
| 90 | 0.029807 | 0 | -0.00745 | 0 | 0.004968 | 0 | -0.00745 | 0 | 0.029807 | 0 |
| 120 | -0.46856 | -0.5 | -0.50582 | -0.5 | -0.4934 | -0.5 | -0.50582 | -0.5 | -0.46856 | -0.5 |
| 150 | -0.83983 | -0.86603 | -0.87709 | -0.86603 | -0.86467 | -0.86603 | -0.87709 | -0.86603 | -0.83983 | -0.86603 |
| 180 | -0.9569 | -1 | -0.99416 | -1 | -0.98174 | -1 | -0.99416 | -1 | -0.9569 | -1 |

Table 18. Accuracy of $y_{\eta\eta} = y_{\eta\eta}(\xi,\eta)$.

| $\eta$ deg | $\xi$= 2.0 Num. | $\xi$= 2.0 Exact | $\xi$= 3.0 Num. | $\xi$= 3.0 Exact | $\xi$= 4.0 Num. | $\xi$= 4.0 Exact | $\xi$= 5.0 Num. | $\xi$= 5.0 Exact | $\xi$= 6.0 Num. | $\xi$= 6.0 Exact |
|---|---|---|---|---|---|---|---|---|---|---|
| 0 | -0.22292 | 0 | -0.12284 | 0 | -0.02276 | 0 | 0.07732 | 0 | 0.177401 | 0 |
| 30 | -0.98479 | -1 | -1.49157 | -1.5 | -1.99836 | -2 | -2.50514 | -2.5 | -3.01192 | -3 |
| 60 | -1.72259 | -1.73205 | -2.5935 | -2.59808 | -3.4644 | -3.4641 | -4.3353 | -4.33013 | -5.2062 | -5.19615 |
| 90 | -2.01531 | -2 | -3.00759 | -3 | -3.99986 | -4 | -4.99213 | -5 | -5.98441 | -6 |
| 120 | -1.72259 | -1.73205 | -2.5935 | -2.59808 | -3.4644 | -3.4641 | -4.3353 | -4.33013 | -5.2062 | -5.19615 |
| 150 | -0.98479 | -1 | -1.49157 | -1.5 | -1.99836 | -2 | -2.50514 | -2.5 | -3.01192 | -3 |
| 180 | -0.22292 | 0 | -0.12284 | 0 | -0.02276 | 0 | 0.07732 | 0 | 0.177401 | 0 |

The numerical of $\phi$ is shown in Fig. 28. The effects of the eccentricity seem reflected properly.

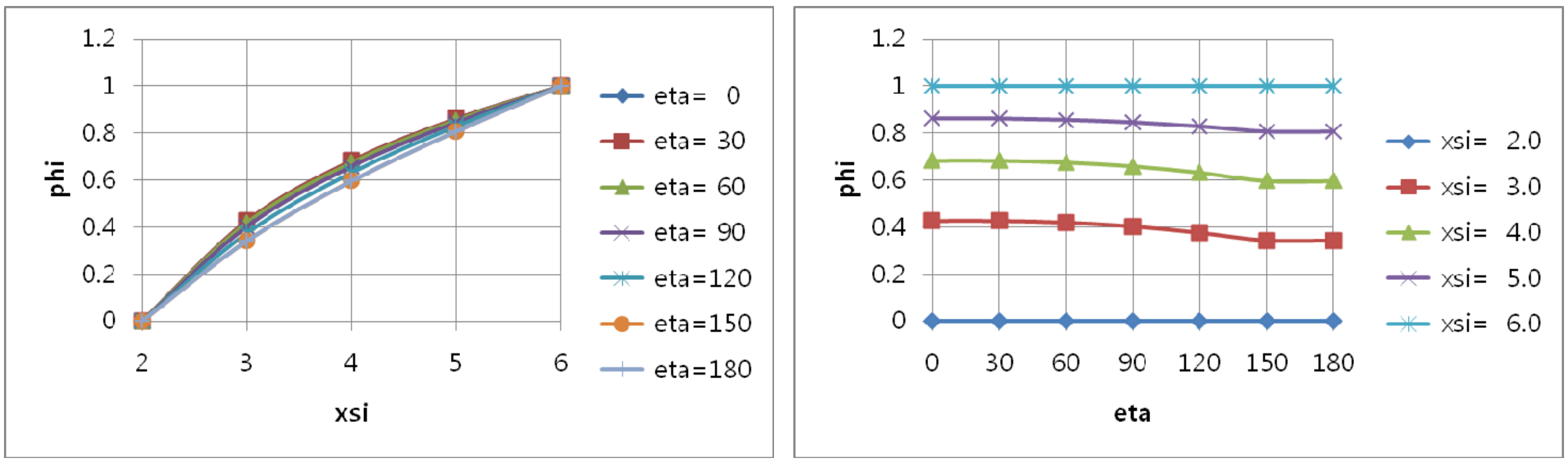

(a) Radial distribution (b) Circumerential distribution

Fig. 28. Distribution of solution $\phi$.

## 6. Conclusions

Graphical introduction of curvilinear coordinates are discussed. For the interpolation of discrete mapping function values, an algebraic method using LSM is used.

Accuracy of mapping is checked for various cases. The generalization ability of the interpolation function is studied. This ability is important to estimate the derivatives of the mapping function. If the parameters are carefully chosen, we can assure the generalization ability.

Not only the mapping of two coordinates but also the application to the solution of boundary value problem is given.

**References**

1. Hiroshi Isshiki, Daisuke Kitazawa, Introduction of Curvilinear Coordinates into Numerical Analysis, Cornell arXiv: 1707.03822, (2017).
2. J. F. Thompson, B. K. Soni, N. P. Weatherill, Handbook of Grid Generation, CRC Press, (1998).
3. C. H. Hu, Slow Drift Damping due to Viscous Force acting on a Cylinder of Floating Structure, Doctor dissertation, Shanghai Jiao Tong University, (1994).